# KERNEL METHODS IN MACHINE LEARNING[1]

By Thomas Hofmann, Bernhard Schölkopf
and Alexander J. Smola

*Darmstadt University of Technology, Max Planck Institute for Biological Cybernetics and National ICT Australia*

We review machine learning methods employing positive definite kernels. These methods formulate learning and estimation problems in a reproducing kernel Hilbert space (RKHS) of functions defined on the data domain, expanded in terms of a kernel. Working in linear spaces of function has the benefit of facilitating the construction and analysis of learning algorithms while at the same time allowing large classes of functions. The latter include nonlinear functions as well as functions defined on nonvectorial data.

We cover a wide range of methods, ranging from binary classifiers to sophisticated methods for estimation with structured data.

**1. Introduction.** Over the last ten years estimation and learning methods utilizing positive definite kernels have become rather popular, particularly in machine learning. Since these methods have a stronger mathematical slant than earlier machine learning methods (e.g., neural networks), there is also significant interest in the statistics and mathematics community for these methods. The present review aims to summarize the state of the art on a conceptual level. In doing so, we build on various sources, including Burges [25], Cristianini and Shawe-Taylor [37], Herbrich [64] and Vapnik [141] and, in particular, Schölkopf and Smola [118], but we also add a fair amount of more recent material which helps unifying the exposition. We have not had space to include proofs; they can be found either in the long version of the present paper (see Hofmann et al. [69]), in the references given or in the above books.

The main idea of all the described methods can be summarized in one paragraph. Traditionally, theory and algorithms of machine learning and

Received December 2005; revised February 2007.
[1]Supported in part by grants of the ARC and by the Pascal Network of Excellence.

statistics has been very well developed for the linear case. Real world data analysis problems, on the other hand, often require nonlinear methods to detect the kind of dependencies that allow successful prediction of properties of interest. By using a positive definite kernel, one can sometimes have the best of both worlds. The kernel corresponds to a dot product in a (usually high-dimensional) feature space. In this space, our estimation methods are linear, but as long as we can formulate everything in terms of kernel evaluations, we never explicitly have to compute in the high-dimensional feature space.

The paper has three main sections: Section 2 deals with fundamental properties of *kernels*, with special emphasis on (conditionally) positive definite kernels and their characterization. We give concrete examples for such kernels and discuss kernels and reproducing kernel Hilbert spaces in the context of regularization. Section 3 presents various approaches for estimating dependencies and analyzing data that make use of kernels. We provide an overview of the problem formulations as well as their solution using convex programming techniques. Finally, Section 4 examines the use of reproducing kernel Hilbert spaces as a means to define statistical models, the focus being on structured, multidimensional responses. We also show how such techniques can be combined with Markov networks as a suitable framework to model dependencies between response variables.

## 2. Kernels.

2.1. *An introductory example.* Suppose we are given empirical data

$$(1) \qquad (x_1, y_1), \ldots, (x_n, y_n) \in \mathcal{X} \times \mathcal{Y}.$$

Here, the domain $\mathcal{X}$ is some nonempty set that the *inputs* (the predictor variables) $x_i$ are taken from; the $y_i \in \mathcal{Y}$ are called *targets* (the response variable). Here and below, $i, j \in [n]$, where we use the notation $[n] := \{1, \ldots, n\}$.

Note that we have not made any assumptions on the domain $\mathcal{X}$ other than it being a set. In order to study the problem of learning, we need additional structure. In learning, we want to be able to *generalize* to unseen data points. In the case of binary pattern recognition, given some new input $x \in \mathcal{X}$, we want to predict the corresponding $y \in \{\pm 1\}$ (more complex output domains $\mathcal{Y}$ will be treated below). Loosely speaking, we want to choose $y$ such that $(x, y)$ is in some sense *similar* to the training examples. To this end, we need similarity measures in $\mathcal{X}$ and in $\{\pm 1\}$. The latter is easier, as two target values can only be identical or different. For the former, we require a function

$$(2) \qquad k : \mathcal{X} \times \mathcal{X} \to \mathbb{R}, \qquad (x, x') \mapsto k(x, x')$$

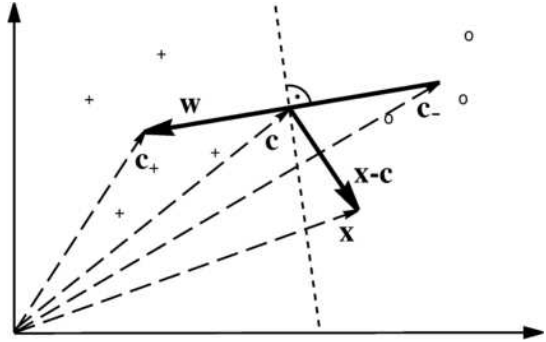

Fig. 1. *A simple geometric classification algorithm: given two classes of points (depicted by "o" and "+"), compute their means $c_+, c_-$ and assign a test input $x$ to the one whose mean is closer. This can be done by looking at the dot product between $x - c$ [where $c = (c_+ + c_-)/2$] and $\mathbf{w} := c_+ - c_-$, which changes sign as the enclosed angle passes through $\pi/2$. Note that the corresponding decision boundary is a hyperplane (the dotted line) orthogonal to $\mathbf{w}$ (from Schölkopf and Smola [118]).*

satisfying, for all $x, x' \in \mathcal{X}$,

$$(3) \qquad k(x, x') = \langle \Phi(x), \Phi(x') \rangle,$$

where $\Phi$ maps into some dot product space $\mathcal{H}$, sometimes called the *feature space*. The similarity measure $k$ is usually called a *kernel*, and $\Phi$ is called its *feature map*.

The advantage of using such a kernel as a similarity measure is that it allows us to construct algorithms in dot product spaces. For instance, consider the following simple classification algorithm, described in Figure 1, where $\mathcal{Y} = \{\pm 1\}$. The idea is to compute the means of the two classes in the feature space, $c_+ = \frac{1}{n_+} \sum_{\{i: y_i = +1\}} \Phi(x_i)$, and $c_- = \frac{1}{n_-} \sum_{\{i: y_i = -1\}} \Phi(x_i)$, where $n_+$ and $n_-$ are the number of examples with positive and negative target values, respectively. We then assign a new point $\Phi(x)$ to the class whose mean is closer to it. This leads to the prediction rule

$$(4) \qquad y = \text{sgn}(\langle \Phi(x), c_+ \rangle - \langle \Phi(x), c_- \rangle + b)$$

with $b = \frac{1}{2}(\|c_-\|^2 - \|c_+\|^2)$. Substituting the expressions for $c_\pm$ yields

$$(5) \quad y = \text{sgn}\left( \frac{1}{n_+} \sum_{\{i: y_i = +1\}} \underbrace{\langle \Phi(x), \Phi(x_i) \rangle}_{k(x, x_i)} - \frac{1}{n_-} \sum_{\{i: y_i = -1\}} \underbrace{\langle \Phi(x), \Phi(x_i) \rangle}_{k(x, x_i)} + b \right),$$

where $b = \frac{1}{2}(\frac{1}{n_-^2} \sum_{\{(i,j): y_i = y_j = -1\}} k(x_i, x_j) - \frac{1}{n_+^2} \sum_{\{(i,j): y_i = y_j = +1\}} k(x_i, x_j))$.

Let us consider one well-known special case of this type of classifier. Assume that the class means have the same distance to the origin (hence, $b = 0$), and that $k(\cdot, x)$ is a density for all $x \in \mathcal{X}$. If the two classes are equally likely and were generated from two probability distributions that are estimated

$$(6) \quad p_+(x) := \frac{1}{n_+} \sum_{\{i:y_i=+1\}} k(x,x_i), \qquad p_-(x) := \frac{1}{n_-} \sum_{\{i:y_i=-1\}} k(x,x_i),$$

then (5) is the estimated Bayes decision rule, plugging in the estimates $p_+$ and $p_-$ for the true densities.

The classifier (5) is closely related to the *Support Vector Machine* (*SVM*) that we will discuss below. It is linear in the feature space (4), while in the input domain, it is represented by a kernel expansion (5). In both cases, the decision boundary is a hyperplane in the feature space; however, the normal vectors [for (4), $w = c_+ - c_-$] are usually rather different.

The normal vector not only characterizes the alignment of the hyperplane, its length can also be used to construct tests for the equality of the two class-generating distributions (Borgwardt et al. [22]).

As an aside, note that if we normalize the targets such that $\hat{y}_i = y_i/|\{j : y_j = y_i\}|$, in which case the $\hat{y}_i$ sum to zero, then $\|w\|^2 = \langle K, \hat{y}\hat{y}^\top \rangle_F$, where $\langle \cdot, \cdot \rangle_F$ is the Frobenius dot product. If the two classes have equal size, then up to a scaling factor involving $\|K\|_2$ and $n$, this equals the *kernel-target alignment* defined by Cristianini et al. [38].

2.2. *Positive definite kernels.* We have required that a kernel satisfy (3), that is, correspond to a dot product in some dot product space. In the present section we show that the class of kernels that can be written in the form (3) coincides with the class of positive definite kernels. This has far-reaching consequences. There are examples of positive definite kernels which can be evaluated efficiently even though they correspond to dot products in infinite dimensional dot product spaces. In such cases, substituting $k(x,x')$ for $\langle \Phi(x), \Phi(x') \rangle$, as we have done in (5), is crucial. In the machine learning community, this substitution is called the *kernel trick*.

DEFINITION 1 (Gram matrix). Given a kernel $k$ and inputs $x_1, \ldots, x_n \in \mathcal{X}$, the $n \times n$ matrix

$$(7) \qquad K := (k(x_i, x_j))_{ij}$$

is called the Gram matrix (or kernel matrix) of $k$ with respect to $x_1, \ldots, x_n$.

DEFINITION 2 (Positive definite matrix). A real $n \times n$ symmetric matrix $K_{ij}$ satisfying

$$(8) \qquad \sum_{i,j} c_i c_j K_{ij} \geq 0$$

for all $c_i \in \mathbb{R}$ is called *positive definite*. If equality in (8) only occurs for $c_1 = \cdots = c_n = 0$, then we shall call the matrix *strictly positive definite*.

DEFINITION 3 (Positive definite kernel). Let $\mathcal{X}$ be a nonempty set. A function $k: \mathcal{X} \times \mathcal{X} \to \mathbb{R}$ which for all $n \in \mathbb{N}, x_i \in \mathcal{X}$, $i \in [n]$ gives rise to a positive definite Gram matrix is called a *positive definite kernel*. A function $k: \mathcal{X} \times \mathcal{X} \to \mathbb{R}$ which for all $n \in \mathbb{N}$ and distinct $x_i \in \mathcal{X}$ gives rise to a strictly positive definite Gram matrix is called a *strictly positive definite kernel*.

Occasionally, we shall refer to positive definite kernels simply as *kernels*. Note that, for simplicity, we have restricted ourselves to the case of real valued kernels. However, with small changes, the below will also hold for the complex valued case.

Since $\sum_{i,j} c_i c_j \langle \Phi(x_i), \Phi(x_j) \rangle = \langle \sum_i c_i \Phi(x_i), \sum_j c_j \Phi(x_j) \rangle \geq 0$, kernels of the form (3) are positive definite for any choice of $\Phi$. In particular, if $\mathcal{X}$ is already a dot product space, we may choose $\Phi$ to be the identity. Kernels can thus be regarded as generalized dot products. While they are not generally bilinear, they share important properties with dot products, such as the Cauchy–Schwarz inequality: If $k$ is a positive definite kernel, and $x_1, x_2 \in \mathcal{X}$, then

$$(9) \qquad k(x_1, x_2)^2 \leq k(x_1, x_1) \cdot k(x_2, x_2).$$

2.2.1. *Construction of the reproducing kernel Hilbert space.* We now define a map from $\mathcal{X}$ into the space of functions mapping $\mathcal{X}$ into $\mathbb{R}$, denoted as $\mathbb{R}^{\mathcal{X}}$, via

$$(10) \qquad \Phi: \mathcal{X} \to \mathbb{R}^{\mathcal{X}} \qquad \text{where } x \mapsto k(\cdot, x).$$

Here, $\Phi(x) = k(\cdot, x)$ denotes the function that assigns the value $k(x', x)$ to $x' \in \mathcal{X}$.

We next construct a dot product space containing the images of the inputs under $\Phi$. To this end, we first turn it into a vector space by forming linear combinations

$$(11) \qquad f(\cdot) = \sum_{i=1}^{n} \alpha_i k(\cdot, x_i).$$

Here, $n \in \mathbb{N}$, $\alpha_i \in \mathbb{R}$ and $x_i \in \mathcal{X}$ are arbitrary.

Next, we define a dot product between $f$ and another function $g(\cdot) = \sum_{j=1}^{n'} \beta_j k(\cdot, x'_j)$ (with $n' \in \mathbb{N}$, $\beta_j \in \mathbb{R}$ and $x'_j \in \mathcal{X}$) as

$$(12) \qquad \langle f, g \rangle := \sum_{i=1}^{n} \sum_{j=1}^{n'} \alpha_i \beta_j k(x_i, x'_j).$$

To see that this is well defined although it contains the expansion coefficients and points, note that $\langle f, g \rangle = \sum_{j=1}^{n'} \beta_j f(x'_j)$. The latter, however, does not depend on the particular expansion of $f$. Similarly, for $g$, note that $\langle f, g \rangle = \sum_{i=1}^{n} \alpha_i g(x_i)$. This also shows that $\langle \cdot, \cdot \rangle$ is bilinear. It is symmetric, as $\langle f, g \rangle =$ $\langle g, f \rangle$. Moreover, it is positive definite, since positive definiteness of $k$ implies that, for any function $f$, written as (11), we have

$$\langle f, f \rangle = \sum_{i,j=1}^{n} \alpha_i \alpha_j k(x_i, x_j) \geq 0. \tag{13}$$

Next, note that given functions $f_1, \ldots, f_p$, and coefficients $\gamma_1, \ldots, \gamma_p \in \mathbb{R}$, we have

$$\sum_{i,j=1}^{p} \gamma_i \gamma_j \langle f_i, f_j \rangle = \left\langle \sum_{i=1}^{p} \gamma_i f_i, \sum_{j=1}^{p} \gamma_j f_j \right\rangle \geq 0. \tag{14}$$

Here, the equality follows from the bilinearity of $\langle \cdot, \cdot \rangle$, and the right-hand inequality from (13).

By (14), $\langle \cdot, \cdot \rangle$ is a positive definite kernel, defined on our vector space of functions. For the last step in proving that it even is a dot product, we note that, by (12), for all functions (11),

$$(15) \quad \langle k(\cdot, x), f \rangle = f(x) \quad \text{and, in particular,} \quad \langle k(\cdot, x), k(\cdot, x') \rangle = k(x, x').$$

By virtue of these properties, $k$ is called a *reproducing kernel* (Aronszajn [7]).

Due to (15) and (9), we have

$$|f(x)|^2 = |\langle k(\cdot, x), f \rangle|^2 \leq k(x, x) \cdot \langle f, f \rangle. \tag{16}$$

By this inequality, $\langle f, f \rangle = 0$ implies $f = 0$, which is the last property that was left to prove in order to establish that $\langle \cdot, \cdot \rangle$ is a dot product.

Skipping some details, we add that one can complete the space of functions (11) in the norm corresponding to the dot product, and thus gets a Hilbert space $\mathcal{H}$, called a *reproducing kernel Hilbert space* (*RKHS*).

One can define a RKHS as a Hilbert space $\mathcal{H}$ of functions on a set $\mathcal{X}$ with the property that, for all $x \in \mathcal{X}$ and $f \in \mathcal{H}$, the point evaluations $f \mapsto f(x)$ are continuous linear functionals [in particular, all point values $f(x)$ are well defined, which already distinguishes RKHSs from many $L_2$ Hilbert spaces]. From the point evaluation functional, one can then construct the reproducing kernel using the Riesz representation theorem. The Moore–Aronszajn theorem (Aronszajn [7]) states that, for every positive definite kernel on $\mathcal{X} \times \mathcal{X}$, there exists a unique RKHS and vice versa.

There is an analogue of the kernel trick for distances rather than dot products, that is, dissimilarities rather than similarities. This leads to the larger class of *conditionally positive definite kernels*. Those kernels are defined just like positive definite ones, with the one difference being that their Gram matrices need to satisfy (8) only subject to

$$\sum_{i=1}^{n} c_i = 0. \tag{17}$$

Interestingly, it turns out that many kernel algorithms, including SVMs and kernel PCA (see Section 3), can be applied also with this larger class of kernels, due to their being translation invariant in feature space (Hein et al. [63] and Schölkopf and Smola [118]).

We conclude this section with a note on terminology. In the early years of kernel machine learning research, it was not the notion of positive definite kernels that was being used. Instead, researchers considered kernels satisfying the conditions of Mercer's theorem (Mercer [99], see, e.g., Cristianini and Shawe-Taylor [37] and Vapnik [141]). However, while all such kernels do satisfy (3), the converse is not true. Since (3) is what we are interested in, positive definite kernels are thus the right class of kernels to consider.

2.2.2. *Properties of positive definite kernels.* We begin with some closure properties of the set of positive definite kernels.

PROPOSITION 4. *Below, $k_1, k_2, \ldots$ are arbitrary positive definite kernels on $\mathcal{X} \times \mathcal{X}$, where $\mathcal{X}$ is a nonempty set:*

(i) *The set of positive definite kernels is a closed convex cone, that is,* (a) *if $\alpha_1, \alpha_2 \geq 0$, then $\alpha_1 k_1 + \alpha_2 k_2$ is positive definite; and* (b) *if $k(x, x') := \lim_{n \to \infty} k_n(x, x')$ exists for all $x, x'$, then $k$ is positive definite.*

(ii) *The pointwise product $k_1 k_2$ is positive definite.*

(iii) *Assume that for $i = 1, 2$, $k_i$ is a positive definite kernel on $\mathcal{X}_i \times \mathcal{X}_i$, where $\mathcal{X}_i$ is a nonempty set. Then the tensor product $k_1 \otimes k_2$ and the direct sum $k_1 \oplus k_2$ are positive definite kernels on $(\mathcal{X}_1 \times \mathcal{X}_2) \times (\mathcal{X}_1 \times \mathcal{X}_2)$.*

The proofs can be found in Berg et al. [18].

It is reassuring that sums and products of positive definite kernels are positive definite. We will now explain that, loosely speaking, there are no other operations that preserve positive definiteness. To this end, let $C$ denote the set of all functions $\psi \colon \mathbb{R} \to \mathbb{R}$ that map positive definite kernels to (conditionally) positive definite kernels (readers who are not interested in the case of conditionally positive definite kernels may ignore the term in parentheses). We define

$C := \{\psi | k \text{ is a p.d. kernel} \Rightarrow \psi(k) \text{ is a (conditionally) p.d. kernel}\},$

$C' = \{\psi | \text{ for any Hilbert space } \mathcal{F},$

$\psi(\langle x, x' \rangle_\mathcal{F}) \text{ is (conditionally) positive definite}\},$

$C'' = \{\psi | \text{ for all } n \in \mathbb{N} \colon K \text{ is a p.d.}$

$n \times n \text{ matrix } \Rightarrow \psi(K) \text{ is (conditionally) p.d.}\},$

where $\psi(K)$ is the $n \times n$ matrix with elements $\psi(K_{ij})$.

PROPOSITION 5.　$C = C' = C''$.

The following proposition follows from a result of FitzGerald et al. [50] for (conditionally) positive definite matrices; by Proposition 5, it also applies for (conditionally) positive definite kernels, and for functions of dot products. We state the latter case.

PROPOSITION 6.　*Let $\psi : \mathbb{R} \to \mathbb{R}$. Then $\psi(\langle x, x' \rangle_{\mathcal{F}})$ is positive definite for any Hilbert space $\mathcal{F}$ if and only if $\psi$ is real entire of the form*

$$\psi(t) = \sum_{n=0}^{\infty} a_n t^n \tag{18}$$

*with $a_n \geq 0$ for $n \geq 0$.*

*Moreover, $\psi(\langle x, x' \rangle_{\mathcal{F}})$ is conditionally positive definite for any Hilbert space $\mathcal{F}$ if and only if $\psi$ is real entire of the form* (18) *with $a_n \geq 0$ for $n \geq 1$.*

There are further properties of $k$ that can be read off the coefficients $a_n$:

- Steinwart [128] showed that if all $a_n$ are strictly positive, then the kernel of Proposition 6 is *universal* on every compact subset $S$ of $\mathbb{R}^d$ in the sense that its RKHS is dense in the space of continuous functions on $S$ in the $\|\cdot\|_\infty$ norm. For support vector machines using universal kernels, he then shows (universal) consistency (Steinwart [129]). Examples of universal kernels are (19) and (20) below.
- In Lemma 11 we will show that the $a_0$ term does not affect an SVM. Hence, we infer that it is actually sufficient for consistency to have $a_n > 0$ for $n \geq 1$.

We conclude the section with an example of a kernel which is positive definite by Proposition 6. To this end, let $\mathcal{X}$ be a dot product space. The power series expansion of $\psi(x) = e^x$ then tells us that

$$k(x, x') = e^{\langle x, x' \rangle / \sigma^2} \tag{19}$$

is positive definite (Haussler [62]). If we further multiply $k$ with the positive definite kernel $f(x)f(x')$, where $f(x) = e^{-\|x\|^2 / 2\sigma^2}$ and $\sigma > 0$, this leads to the positive definiteness of the Gaussian kernel

$$k'(x, x') = k(x, x') f(x) f(x') = e^{-\|x - x'\|^2 / (2\sigma^2)}. \tag{20}$$

2.2.3. *Properties of positive definite functions.* We now let $\mathcal{X} = \mathbb{R}^d$ and consider positive definite kernels of the form

$$k(x, x') = h(x - x'), \tag{21}$$

in which case $h$ is called a *positive definite function*. The following characterization is due to Bochner [21]. We state it in the form given by Wendland [152].

THEOREM 7. *A continuous function $h$ on $\mathbb{R}^d$ is positive definite if and only if there exists a finite nonnegative Borel measure $\mu$ on $\mathbb{R}^d$ such that*

$$h(x) = \int_{\mathbb{R}^d} e^{-i\langle x, \omega \rangle} \, d\mu(\omega). \tag{22}$$

While normally formulated for complex valued functions, the theorem also holds true for real functions. Note, however, that if we start with an arbitrary nonnegative Borel measure, its Fourier transform may not be real. Real-valued positive definite functions are distinguished by the fact that the corresponding measures $\mu$ are symmetric.

We may normalize $h$ such that $h(0) = 1$ [hence, by (9), $|h(x)| \leq 1$], in which case $\mu$ is a probability measure and $h$ is its characteristic function. For instance, if $\mu$ is a normal distribution of the form $(2\pi/\sigma^2)^{-d/2} e^{-\sigma^2 \|\omega\|^2 / 2} \, d\omega$, then the corresponding positive definite function is the Gaussian $e^{-\|x\|^2/(2\sigma^2)}$; see (20).

Bochner's theorem allows us to interpret the similarity measure $k(x, x') = h(x - x')$ in the frequency domain. The choice of the measure $\mu$ determines which frequency components occur in the kernel. Since the solutions of kernel algorithms will turn out to be finite kernel expansions, the measure $\mu$ will thus determine which frequencies occur in the estimates, that is, it will determine their regularization properties—more on that in Section 2.3.2 below.

Bochner's theorem generalizes earlier work of Mathias, and has itself been generalized in various ways, that is, by Schoenberg [115]. An important generalization considers Abelian semigroups (Berg et al. [18]). In that case, the theorem provides an integral representation of positive definite functions in terms of the semigroup's semicharacters. Further generalizations were given by Krein, for the cases of positive definite kernels and functions with a limited number of negative squares. See Stewart [130] for further details and references.

As above, there are conditions that ensure that the positive definiteness becomes strict.

PROPOSITION 8 (Wendland [152]). *A positive definite function is strictly positive definite if the carrier of the measure in its representation* (22) *contains an open subset.*

This implies that the Gaussian kernel is strictly positive definite.

An important special case of positive definite functions, which includes the Gaussian, are *radial basis functions*. These are functions that can be written as $h(x) = g(\|x\|_2)$ for some function $g : [0, \infty[ \to \mathbb{R}$. They have the property of being invariant under the Euclidean group.

2.2.4. *Examples of kernels.* We have already seen several instances of positive definite kernels, and now intend to complete our selection with a few more examples. In particular, we discuss polynomial kernels, convolution kernels, ANOVA expansions and kernels on documents.

*Polynomial kernels.* From Proposition 4 it is clear that homogeneous polynomial kernels $k(x, x') = \langle x, x' \rangle^p$ are positive definite for $p \in \mathbb{N}$ and $x, x' \in \mathbb{R}^d$. By direct calculation, we can derive the corresponding feature map (Poggio [108]):

$$
\begin{aligned}
\langle x, x' \rangle^p &= \left( \sum_{j=1}^d [x]_j [x']_j \right)^p \\
&= \sum_{j \in [d]^p} [x]_{j_1} \cdots [x]_{j_p} \cdot [x']_{j_1} \cdots [x']_{j_p} = \langle C_p(x), C_p(x') \rangle,
\end{aligned}
\tag{23}
$$

where $C_p$ maps $x \in \mathbb{R}^d$ to the vector $C_p(x)$ whose entries are all possible $p$th degree ordered products of the entries of $x$ (note that $[d]$ is used as a shorthand for $\{1, \ldots, d\}$). The polynomial kernel of degree $p$ thus computes a dot product in the space spanned by all monomials of degree $p$ in the input coordinates. Other useful kernels include the inhomogeneous polynomial,

$$
k(x, x') = (\langle x, x' \rangle + c)^p \qquad \text{where } p \in \mathbb{N} \text{ and } c \geq 0,
\tag{24}
$$

which computes all monomials up to degree $p$.

*Spline kernels.* It is possible to obtain spline functions as a result of kernel expansions (Vapnik et al. [144]) simply by noting that convolution of an even number of indicator functions yields a positive kernel function. Denote by $I_X$ the indicator (or characteristic) function on the set $X$, and denote by $\otimes$ the convolution operation, $(f \otimes g)(x) := \int_{\mathbb{R}^d} f(x') g(x' - x) \, dx'$. Then the B-spline kernels are given by

$$
k(x, x') = B_{2p+1}(x - x') \qquad \text{where } p \in \mathbb{N} \text{ with } B_{i+1} := B_i \otimes B_0.
\tag{25}
$$

Here $B_0$ is the characteristic function on the unit ball in $\mathbb{R}^d$. From the definition of (25), it is obvious that, for odd $m$, we may write $B_m$ as the inner product between functions $B_{m/2}$. Moreover, note that, for even $m$, $B_m$ is not a kernel.

*Convolutions and structures.* Let us now move to kernels defined on structured objects (Haussler [62] and Watkins [151]). Suppose the object $x \in \mathcal{X}$ is composed of $x_p \in \mathcal{X}_p$, where $p \in [P]$ (note that the sets $\mathcal{X}_p$ need not be equal). For instance, consider the string $x = ATG$ and $P = 2$. It is composed of the parts $x_1 = AT$ and $x_2 = G$, or alternatively, of $x_1 = A$ and $x_2 = TG$. Mathematically speaking, the set of "allowed" decompositions can be thought of as a *relation* $R(x_1, \ldots, x_P, x)$, to be read as "$x_1, \ldots, x_P$ constitute the composite object $x$."

Haussler [62] investigated how to define a kernel between composite objects by building on similarity measures that assess their respective *parts*; in other words, kernels $k_p$ defined on $\mathcal{X}_p \times \mathcal{X}_p$. Define the *R-convolution* of $k_1, \ldots, k_P$ as

$$(26) \qquad [k_1 \star \cdots \star k_P](x, x') := \sum_{\bar{x} \in R(x), \bar{x}' \in R(x')} \prod_{p=1}^{P} k_p(\bar{x}_p, \bar{x}'_p),$$

where the sum runs over all possible ways $R(x)$ and $R(x')$ in which we can decompose $x$ into $\bar{x}_1, \ldots, \bar{x}_P$ and $x'$ analogously [here we used the convention that an empty sum equals zero, hence, if either $x$ or $x'$ cannot be decomposed, then $(k_1 \star \cdots \star k_P)(x, x') = 0$]. If there is only a finite number of ways, the relation $R$ is called finite. In this case, it can be shown that the $R$-convolution is a valid kernel (Haussler [62]).

*ANOVA kernels.* Specific examples of convolution kernels are Gaussians and ANOVA kernels (Vapnik [141] and Wahba [148]). To construct an ANOVA kernel, we consider $\mathcal{X} = S^N$ for some set $S$, and kernels $k^{(i)}$ on $S \times S$, where $i = 1, \ldots, N$. For $P = 1, \ldots, N$, the *ANOVA kernel of order $P$* is defined as

$$(27) \qquad k_P(x, x') := \sum_{1 \leq i_1 < \cdots < i_P \leq N} \prod_{p=1}^{P} k^{(i_p)}(x_{i_p}, x'_{i_p}).$$

Note that if $P = N$, the sum consists only of the term for which $(i_1, \ldots, i_P) = (1, \ldots, N)$, and $k$ equals the tensor product $k^{(1)} \otimes \cdots \otimes k^{(N)}$. At the other extreme, if $P = 1$, then the products collapse to one factor each, and $k$ equals the direct sum $k^{(1)} \oplus \cdots \oplus k^{(N)}$. For intermediate values of $P$, we get kernels that lie in between tensor products and direct sums.

ANOVA kernels typically use some moderate value of $P$, which specifies the order of the interactions between attributes $x_{i_p}$ that we are interested in. The sum then runs over the numerous terms that take into account interactions of order $P$; fortunately, the computational cost can be reduced to $O(Pd)$ cost by utilizing recurrent procedures for the kernel evaluation. ANOVA kernels have been shown to work rather well in multi-dimensional SV regression problems (Stitson et al. [131]).

*Bag of words.* One way in which SVMs have been used for text categorization (Joachims [77]) is the *bag-of-words* representation. This maps a given text to a sparse vector, where each component corresponds to a word, and a component is set to one (or some other number) whenever the related word occurs in the text. Using an efficient sparse representation, the dot product between two such vectors can be computed quickly. Furthermore, this dot product is by construction a valid kernel, referred to as a *sparse vector kernel*. One of its shortcomings, however, is that it does not take into account the word ordering of a document. Other sparse vector kernels are also conceivable, such as one that maps a text to the set of pairs of words that are in the same sentence (Joachims [77] and Watkins [151]).

*n-grams and suffix trees.* A more sophisticated way of dealing with string data was proposed by Haussler [62] and Watkins [151]. The basic idea is as described above for general structured objects (26): Compare the strings by means of the substrings they contain. The more substrings two strings have in common, the more similar they are. The substrings need not always be contiguous; that said, the further apart the first and last element of a substring are, the less weight should be given to the similarity. Depending on the specific choice of a similarity measure, it is possible to define more or less efficient kernels which compute the dot product in the feature space spanned by *all* substrings of documents.

Consider a finite alphabet $\Sigma$, the set of all strings of length $n$, $\Sigma^n$, and the set of all finite strings, $\Sigma^* := \bigcup_{n=0}^{\infty} \Sigma^n$. The length of a string $s \in \Sigma^*$ is denoted by $|s|$, and its elements by $s(1)\ldots s(|s|)$; the concatenation of $s$ and $t \in \Sigma^*$ is written $st$. Denote by

$$k(x, x') = \sum_s \#(x, s) \#(x', s) c_s$$

a string kernel computed from exact matches. Here $\#(x, s)$ is the number of occurrences of $s$ in $x$ and $c_s \geq 0$.

Vishwanathan and Smola [146] provide an algorithm using suffix trees, which allows one to compute for arbitrary $c_s$ the value of the kernel $k(x, x')$ in $O(|x| + |x'|)$ time and memory. Moreover, also $f(x) = \langle w, \Phi(x) \rangle$ can be computed in $O(|x|)$ time if preprocessing linear in the size of the support vectors is carried out. These kernels are then applied to function prediction (according to the gene ontology) of proteins using only their sequence information. Another prominent application of string kernels is in the field of splice form prediction and gene finding (Rätsch et al. [112]).

For inexact matches of a limited degree, typically up to $\epsilon = 3$, and strings of bounded length, a similar data structure can be built by explicitly generating a dictionary of strings and their neighborhood in terms of a Hamming distance (Leslie et al. [92]). These kernels are defined by replacing $\#(x, s)$ by a mismatch function $\#(x, s, \epsilon)$ which reports the number of approximate occurrences of $s$ in $x$. By trading off computational complexity with storage (hence, the restriction to small numbers of mismatches), essentially linear-time algorithms can be designed. Whether a general purpose algorithm exists which allows for efficient comparisons of strings with mismatches in linear time is still an open question.

*Mismatch kernels.* In the general case it is only possible to find algorithms whose complexity is linear in the lengths of the documents being compared, and the length of the substrings, that is, $O(|x| \cdot |x'|)$ or worse. We now describe such a kernel with a specific choice of weights (Cristianini and Shawe-Taylor [37] and Watkins [151]).

Let us now form subsequences $u$ of strings. Given an index sequence $\mathbf{i} := (i_1, \ldots, i_{|u|})$ with $1 \leq i_1 < \cdots < i_{|u|} \leq |s|$, we define $u := s(\mathbf{i}) := s(i_1) \ldots s(i_{|u|})$. We call $l(\mathbf{i}) := i_{|u|} - i_1 + 1$ *the length of the subsequence in $s$*. Note that if $\mathbf{i}$ is not contiguous, then $l(\mathbf{i}) > |u|$.

The feature space built from strings of length $n$ is defined to be $\mathcal{H}_n := \mathbb{R}^{(\Sigma^n)}$. This notation means that the space has one dimension (or coordinate) for each element of $\Sigma^n$, labeled by that element (equivalently, we can think of it as the space of all real-valued functions on $\Sigma^n$). We can thus describe the feature map coordinate-wise for each $u \in \Sigma^n$ via

$$(28) \qquad [\Phi_n(s)]_u := \sum_{\mathbf{i}:s(\mathbf{i})=u} \lambda^{l(\mathbf{i})}.$$

Here, $0 < \lambda \leq 1$ is a decay parameter: The larger the length of the subsequence in $s$, the smaller the respective contribution to $[\Phi_n(s)]_u$. The sum runs over all subsequences of $s$ which equal $u$.

For instance, consider a dimension of $\mathcal{H}_3$ spanned (i.e., labeled) by the string asd. In this case we have $[\Phi_3(\text{N}\underline{\text{asd}}\text{aq})]_{\text{asd}} = \lambda^3$, while $[\Phi_3(\text{lass das})]_{\text{asd}} = 2\lambda^5$. In the first string, asd is a contiguous substring. In the second string, it appears twice as a noncontiguous substring of length 5 in lass das, the two occurrences are l<u>ass</u> <u>d</u>as and l<u>as</u>s <u>da</u>s.

The kernel induced by the map $\Phi_n$ takes the form

$$(29) \qquad k_n(s,t) = \sum_{u \in \Sigma^n} [\Phi_n(s)]_u [\Phi_n(t)]_u = \sum_{u \in \Sigma^n} \sum_{(\mathbf{i},\mathbf{j}):s(\mathbf{i})=t(\mathbf{j})=u} \lambda^{l(\mathbf{i})} \lambda^{l(\mathbf{j})}.$$

The string kernel $k_n$ can be computed using dynamic programming; see Watkins [151].

The above kernels on string, suffix-tree, mismatch and tree kernels have been used in sequence analysis. This includes applications in document analysis and categorization, spam filtering, function prediction in proteins, annotations of dna sequences for the detection of introns and exons, named entity tagging of documents and the construction of parse trees.

*Locality improved kernels.* It is possible to adjust kernels to the structure of spatial data. Recall the Gaussian RBF and polynomial kernels. When applied to an image, it makes no difference whether one uses as $x$ the image or a version of $x$ where all locations of the pixels have been permuted. This indicates that function space on $\mathcal{X}$ induced by $k$ does not take advantage of the *locality* properties of the data.

By taking advantage of the local structure, estimates can be improved. On biological sequences (Zien et al. [157]) one may assign more weight to the entries of the sequence close to the location where estimates should occur.

For images, local interactions between image patches need to be considered. One way is to use the *pyramidal* kernel (DeCoste and Schölkopf [44] and Schölkopf [116]). It takes inner products between corresponding image patches, then raises the latter to some power $p_1$, and finally raises their sum to another power $p_2$. While the overall degree of this kernel is $p_1 p_2$, the first factor $p_1$ only captures short range interactions.

*Tree kernels.* We now discuss similarity measures on more structured objects. For trees Collins and Duffy [31] propose a decomposition method which maps a tree $x$ into its set of subtrees. The kernel between two trees $x, x'$ is then computed by taking a weighted sum of all terms between both trees. In particular, Collins and Duffy [31] show a quadratic time algorithm, that is, $O(|x| \cdot |x'|)$ to compute this expression, where $|x|$ is the number of nodes of the tree. When restricting the sum to all proper rooted subtrees, it is possible to reduce the computational cost to $O(|x| + |x'|)$ time by means of a tree to string conversion (Vishwanathan and Smola [146]).

*Graph kernels.* Graphs pose a twofold challenge: one may both design a kernel *on* vertices of them and also a kernel *between* them. In the former case, the graph itself becomes the object defining the metric between the vertices. See Gärtner [56] and Kashima et al. [82] for details on the latter. In the following we discuss kernels *on* graphs.

Denote by $W \in \mathbb{R}^{n \times n}$ the adjacency matrix of a graph with $W_{ij} > 0$ if an edge between $i, j$ exists. Moreover, assume for simplicity that the graph is undirected, that is, $W^\top = W$. Denote by $L = D - W$ the graph Laplacian and by $\tilde{L} = \mathbf{1} - D^{-1/2} W D^{-1/2}$ the normalized graph Laplacian. Here $D$ is a diagonal matrix with $D_{ii} = \sum_j W_{ij}$ denoting the degree of vertex $i$.

Fiedler [49] showed that, the second largest eigenvector of $L$ approximately decomposes the graph into two parts according to their sign. The other large eigenvectors partition the graph into correspondingly smaller portions. $L$ arises from the fact that for a function $f$ defined on the vertices of the graph $\sum_{i,j} (f(i) - f(j))^2 = 2 f^\top L f$.

Finally, Smola and Kondor [125] show that, under mild conditions and up to rescaling, $L$ is the only quadratic permutation invariant form which can be obtained as a linear function of $W$.

Hence, it is reasonable to consider kernel matrices $K$ obtained from $L$ (and $\tilde{L}$). Smola and Kondor [125] suggest kernels $K = r(L)$ or $K = r(\tilde{L})$, which have desirable smoothness properties. Here $r:[0,\infty) \to [0,\infty)$ is a monotonically decreasing function. Popular choices include

$$(30) \quad r(\xi) = \exp(-\lambda \xi) \quad \text{diffusion kernel,}$$

$$(31) \quad r(\xi) = (\xi + \lambda)^{-1} \quad \text{regularized graph Laplacian,}$$

$$(32) \quad r(\xi) = (\lambda - \xi)^p \quad p\text{-step random walk,}$$

where $\lambda > 0$ is chosen such as to reflect the amount of diffusion in (30), the degree of regularization in (31) or the weighting of steps within a random walk (32) respectively. Equation (30) was proposed by Kondor and Lafferty [87]. In Section 2.3.2 we will discuss the connection between regularization operators and kernels in $\mathbb{R}^n$. Without going into details, the function $r(\xi)$ describes the smoothness properties on the graph and $L$ plays the role of the Laplace operator.

*Kernels on sets and subspaces.* Whenever each observation $x_i$ consists of a *set* of instances, we may use a range of methods to capture the specific properties of these sets (for an overview, see Vishwanathan et al. [147]):

- Take the average of the elements of the set in feature space, that is, $\phi(x_i) = \frac{1}{n}\sum_j \phi(x_{ij})$. This yields good performance in the area of multi-instance learning.
- Jebara and Kondor [75] extend the idea by dealing with distributions $p_i(x)$ such that $\phi(x_i) = \mathbf{E}[\phi(x)]$, where $x \sim p_i(x)$. They apply it to image classification with missing pixels.
- Alternatively, one can study angles enclosed by *subspaces* spanned by the observations. In a nutshell, if $U, U'$ denote the orthogonal matrices spanning the subspaces of $x$ and $x'$ respectively, then $k(x,x') = \det U^\top U'$.

*Fisher kernels.* [74] have designed kernels building on probability density models $p(x|\theta)$. Denote by

$$(33) \quad U_\theta(x) := -\partial_\theta \log p(x|\theta),$$

$$(34) \quad I := \mathbf{E}_x[U_\theta(x) U_\theta^\top(x)],$$

the Fisher scores and the Fisher information matrix respectively. Note that for maximum likelihood estimators $\mathbf{E}_x[U_\theta(x)] = 0$ and, therefore, $I$ is the covariance of $U_\theta(x)$. The Fisher kernel is defined as

$$(35) \quad k(x,x') := U_\theta^\top(x) I^{-1} U_\theta(x') \text{ or } k(x,x') := U_\theta^\top(x) U_\theta(x')$$

depending on whether we study the normalized or the unnormalized kernel respectively.

In addition to that, it has several attractive theoretical properties: Oliver et al. [104] show that estimation using the normalized Fisher kernel corresponds to estimation subject to a regularization on the $L_2(p(\cdot|\theta))$ norm.

Moreover, in the context of exponential families (see Section 4.1 for a more detailed discussion) where $p(x|\theta) = \exp(\langle \phi(x), \theta \rangle - g(\theta))$, we have

$$(36) \qquad k(x, x') = [\phi(x) - \partial_\theta g(\theta)][\phi(x') - \partial_\theta g(\theta)]$$

for the unnormalized Fisher kernel. This means that up to centering by $\partial_\theta g(\theta)$ the Fisher kernel is identical to the kernel arising from the inner product of the sufficient statistics $\phi(x)$. This is not a coincidence. In fact, in our analysis of nonparametric exponential families we will encounter this fact several times (cf. Section 4 for further details). Moreover, note that the centering is immaterial, as can be seen in Lemma 11.

The above overview of kernel design is by no means complete. The reader is referred to books of Bakir et al. [9], Cristianini and Shawe-Taylor [37], Herbrich [64], Joachims [77], Schölkopf and Smola [118], Schölkopf [121] and Shawe-Taylor and Cristianini [123] for further examples and details.

### 2.3. Kernel function classes.

#### 2.3.1. The representer theorem.

From kernels, we now move to functions that can be expressed in terms of kernel expansions. The representer theorem (Kimeldorf and Wahba [85] and Schölkopf and Smola [118]) shows that solutions of a large class of optimization problems can be expressed as kernel expansions over the sample points. As above, $\mathcal{H}$ is the RKHS associated to the kernel $k$.

THEOREM 9 (Representer theorem). *Denote by $\Omega:[0,\infty) \to \mathbb{R}$ a strictly monotonic increasing function, by $\mathcal{X}$ a set, and by $c:(\mathcal{X} \times \mathbb{R}^2)^n \to \mathbb{R} \cup \{\infty\}$ an arbitrary loss function. Then each minimizer $f \in \mathcal{H}$ of the regularized risk functional*

$$(37) \qquad c((x_1, y_1, f(x_1)), \ldots, (x_n, y_n, f(x_n))) + \Omega(\|f\|_\mathcal{H}^2)$$

*admits a representation of the form*

$$(38) \qquad f(x) = \sum_{i=1}^n \alpha_i k(x_i, x).$$

Monotonicity of $\Omega$ does not prevent the regularized risk functional (37) from having multiple local minima. To ensure a global minimum, we would need to require convexity. If we discard the strictness of the monotonicity, then it no longer follows that each minimizer of the regularized risk admits an expansion (38); it still follows, however, that there is always another solution that is as good, and that *does* admit the expansion.

The significance of the representer theorem is that although we might be trying to solve an optimization problem in an infinite-dimensional space $\mathcal{H}$, containing linear combinations of kernels centered on *arbitrary* points of $\mathcal{X}$, it states that the solution lies in the span of $n$ particular kernels—those centered on the training points. We will encounter (38) again further below, where it is called the *Support Vector expansion*. For suitable choices of loss functions, many of the $\alpha_i$ often equal 0.

Despite the finiteness of the representation in (38), it can often be the case that the number of terms in the expansion is too large in practice. This can be problematic in practice, since the time required to evaluate (38) is proportional to the number of terms. One can reduce this number by computing a reduced representation which approximates the original one in the RKHS norm (e.g., Schölkopf and Smola [118]).

2.3.2. *Regularization properties.* The regularizer $\|f\|_{\mathcal{H}}^2$ used in Theorem 9, which is what distinguishes SVMs from many other regularized function estimators (e.g., based on coefficient based $L_1$ regularizers, such as the Lasso (Tibshirani [135]) or linear programming machines (Schölkopf and Smola [118])), stems from the dot product $\langle f, f \rangle_k$ in the RKHS $\mathcal{H}$ associated with a positive definite kernel. The nature and implications of this regularizer, however, are not obvious and we shall now provide an analysis in the Fourier domain. It turns out that if the kernel is translation invariant, then its Fourier transform allows us to characterize how the different frequency components of $f$ contribute to the value of $\|f\|_{\mathcal{H}}^2$. Our exposition will be informal (see also Poggio and Girosi [109] and Smola et al. [127]), and we will implicitly assume that all integrals are over $\mathbb{R}^d$ and exist, and that the operators are well defined.

We will rewrite the RKHS dot product as

$$\langle f, g \rangle_k = \langle \Upsilon f, \Upsilon g \rangle = \langle \Upsilon^2 f, g \rangle, \tag{39}$$

where $\Upsilon$ is a positive (and thus symmetric) operator mapping $\mathcal{H}$ into a function space endowed with the usual dot product

$$\langle f, g \rangle = \int f(x) \overline{g(x)} \, dx. \tag{40}$$

Rather than (39), we consider the equivalent condition (cf. Section 2.2.1)

$$\langle k(x, \cdot), k(x', \cdot) \rangle_k = \langle \Upsilon k(x, \cdot), \Upsilon k(x', \cdot) \rangle = \langle \Upsilon^2 k(x, \cdot), k(x', \cdot) \rangle. \tag{41}$$

If $k(x, \cdot)$ is a *Green function* of $\Upsilon^2$, we have

$$\langle \Upsilon^2 k(x, \cdot), k(x', \cdot) \rangle = \langle \delta_x, k(x', \cdot) \rangle = k(x, x'), \tag{42}$$

which by the reproducing property (15) amounts to the desired equality (41).

For conditionally positive definite kernels, a similar correspondence can be established, with a regularization operator whose null space is spanned by a set of functions which are not regularized [in the case (17), which is sometimes called *conditionally positive definite of order* 1, these are the constants].

We now consider the particular case where the kernel can be written $k(x, x') = h(x - x')$ with a continuous strictly positive definite function $h \in L_1(\mathbb{R}^d)$ (cf. Section 2.2.3). A variation of Bochner's theorem, stated by Wendland [152], then tells us that the measure corresponding to $h$ has a nonvanishing density $\upsilon$ with respect to the Lebesgue measure, that is, that $k$ can be written as

$$(43) \qquad k(x, x') = \int e^{-i\langle x-x', \omega \rangle} \upsilon(\omega) \, d\omega = \int e^{-i\langle x, \omega \rangle} \overline{e^{-i\langle x', \omega \rangle}} \upsilon(\omega) \, d\omega.$$

We would like to rewrite this as $\langle \Upsilon k(x, \cdot), \Upsilon k(x', \cdot) \rangle$ for some linear operator $\Upsilon$. It turns out that a multiplication operator in the Fourier domain will do the job. To this end, recall the $d$-dimensional Fourier transform, given by

$$(44) \qquad F[f](\omega) := (2\pi)^{-d/2} \int f(x) e^{-i\langle x, \omega \rangle} \, dx,$$

$$(45) \qquad \text{with the inverse } F^{-1}[f](x) = (2\pi)^{-d/2} \int f(\omega) e^{i\langle x, \omega \rangle} \, d\omega.$$

Next, compute the Fourier transform of $k$ as

$$(46) \quad \begin{aligned} F[k(x, \cdot)](\omega) &= (2\pi)^{-d/2} \int \int (\upsilon(\omega') e^{-i\langle x, \omega' \rangle}) e^{i\langle x', \omega' \rangle} \, d\omega' e^{-i\langle x', \omega \rangle} \, dx' \\ &= (2\pi)^{d/2} \upsilon(\omega) e^{-i\langle x, \omega \rangle}. \end{aligned}$$

Hence, we can rewrite (43) as

$$(47) \qquad k(x, x') = (2\pi)^{-d} \int \frac{F[k(x, \cdot)](\omega) \overline{F[k(x', \cdot)](\omega)}}{\upsilon(\omega)} \, d\omega.$$

If our regularization operator maps

$$(48) \qquad \Upsilon: f \mapsto (2\pi)^{-d/2} \upsilon^{-1/2} F[f],$$

we thus have

$$(49) \qquad k(x, x') = \int (\Upsilon k(x, \cdot))(\omega) \overline{(\Upsilon k(x', \cdot))(\omega)} \, d\omega,$$

that is, our desired identity (41) holds true.

As required in (39), we can thus interpret the dot product $\langle f, g \rangle_k$ in the RKHS as a dot product $\int (\Upsilon f)(\omega) \overline{(\Upsilon g)(\omega)} \, d\omega$. This allows us to understand regularization properties of $k$ in terms of its (scaled) Fourier transform $\upsilon(\omega)$. *Small* values of $\upsilon(\omega)$ amplify the corresponding frequencies in (48). Penalizing $\langle f, f \rangle_k$ thus amounts to a *strong* attenuation of the corresponding frequencies. Hence, small values of $\upsilon(\omega)$ for large $\|\omega\|$ are desirable, since high-frequency components of $F[f]$ correspond to rapid changes in $f$. It follows that $\upsilon(\omega)$ describes the filter properties of the corresponding regularization operator $\Upsilon$. In view of our comments following Theorem 7, we can translate this insight into probabilistic terms: if the probability measure $\frac{\upsilon(\omega)\,d\omega}{\int \upsilon(\omega)\,d\omega}$ describes the desired filter properties, then the natural translation invariant kernel to use is the characteristic function of the measure.

2.3.3. *Remarks and notes.* The notion of kernels as dot products in Hilbert spaces was brought to the field of machine learning by Aizerman et al. [1], Boser at al. [23], Schölkopf at al. [119] and Vapnik [141]. Aizerman et al. [1] used kernels as a tool in a convergence proof, allowing them to apply the Perceptron convergence theorem to their class of potential function algorithms. To the best of our knowledge, Boser et al. [23] were the first to use kernels to construct a nonlinear estimation algorithm, the hard margin predecessor of the Support Vector Machine, from its linear counterpart, the *generalized portrait* (Vapnik [139] and Vapnik and Lerner [145]). While all these uses were limited to kernels defined on vectorial data, Schölkopf [116] observed that this restriction is unnecessary, and nontrivial kernels on other data types were proposed by Haussler [62] and Watkins [151]. Schölkopf et al. [119] applied the kernel trick to generalize principal component analysis and pointed out the (in retrospect obvious) fact that any algorithm which only uses the data via dot products can be generalized using kernels.

In addition to the above uses of positive definite kernels in machine learning, there has been a parallel, and partly earlier development in the field of statistics, where such kernels have been used, for instance, for time series analysis (Parzen [106]), as well as regression estimation and the solution of inverse problems (Wahba [148]).

In probability theory, positive definite kernels have also been studied in depth since they arise as covariance kernels of stochastic processes; see, for example, Loève [93]. This connection is heavily being used in a subset of the machine learning community interested in prediction with Gaussian processes (Rasmussen and Williams [111]).

In functional analysis, the problem of Hilbert space representations of kernels has been studied in great detail; a good reference is Berg at al. [18]; indeed, a large part of the material in the present section is based on that work. Interestingly, it seems that for a fairly long time, there have been two separate strands of development (Stewart [130]). One of them was the study of positive definite functions, which started later but seems to have been unaware of the fact that it considered a special case of positive definite kernels. The latter was initiated by Hilbert [67] and Mercer [99], and was pursued, for instance, by Schoenberg [115]. Hilbert calls a kernel $k$ *definit* if

$$(50) \qquad \int_a^b \int_a^b k(x, x') f(x) f(x') \, dx \, dx' > 0$$

for all nonzero continuous functions $f$, and shows that all eigenvalues of the corresponding integral operator $f \mapsto \int_a^b k(x, \cdot) f(x) \, dx$ are then positive. If $k$ satisfies the condition (50) subject to the constraint that $\int_a^b f(x) g(x) \, dx = 0$, for some fixed function $g$, Hilbert calls it *relativ definit*. For that case, he shows that $k$ has at most one negative eigenvalue. Note that if $f$ is chosen to be constant, then this notion is closely related to the one of conditionally positive definite kernels; see (17). For further historical details, see the review of Stewart [130] or Berg at al. [18].

**3. Convex programming methods for estimation.** As we saw, kernels can be used both for the purpose of describing nonlinear functions subject to smoothness constraints and for the purpose of computing inner products in some feature space efficiently. In this section we focus on the latter and how it allows us to design methods of estimation based on the geometry of the problems at hand.

Unless stated otherwise, $\mathbf{E}[\cdot]$ denotes the expectation with respect to all random variables of the argument. Subscripts, such as $\mathbf{E}_X[\cdot]$, indicate that the expectation is taken over $X$. We will omit them wherever obvious. Finally, we will refer to $\mathbf{E}_{\text{emp}}[\cdot]$ as the empirical average with respect to an $n$-sample. Given a sample $\mathcal{S} := \{(x_1, y_1), \ldots, (x_n, y_n)\} \subseteq \mathcal{X} \times \mathcal{Y}$, we now aim at finding an affine function $f(x) = \langle w, \phi(x) \rangle + b$ or in some cases a function $f(x, y) = \langle \phi(x, y), w \rangle$ such that the empirical risk on $\mathcal{S}$ is minimized. In the binary classification case this means that we want to maximize the agreement between $\text{sgn } f(x)$ and $y$.

- Minimization of the empirical risk with respect to $(w, b)$ is NP-hard (Minsky and Papert [101]). In fact, Ben-David et al. [15] show that even approximately minimizing the empirical risk is NP-hard, not only for linear function classes but also for spheres and other simple geometrical objects. This means that even if the statistical challenges could be solved, we still would be confronted with a formidable algorithmic problem.
- The indicator function $\{y f(x) < 0\}$ is discontinuous and even small changes in $f$ may lead to large changes in both empirical and expected risk. Properties of such functions can be captured by the VC-dimension (Vapnik and Chervonenkis [142]), that is, the maximum number of observations which can be labeled in an arbitrary fashion by functions of the class. Necessary and sufficient conditions for estimation can be stated in these

KERNEL METHODS IN MACHINE LEARNING 21

terms (Vapnik and Chervonenkis [143]). However, much tighter bounds can be obtained by also using the scale of the class (Alon et al. [3]). In fact, there exist function classes parameterized by a single scalar which have infinite VC-dimension (Vapnik [140]).

Given the difficulty arising from minimizing the empirical risk, we now discuss algorithms which minimize an upper bound on the empirical risk, while providing good computational properties and consistency of the estimators. A discussion of the statistical properties follows in Section 3.6.

3.1. *Support vector classification.* Assume that $\mathcal{S}$ is linearly separable, that is, there exists a linear function $f(x)$ such that $\operatorname{sgn} yf(x) = 1$ on $\mathcal{S}$. In this case, the task of finding a large margin separating hyperplane can be viewed as one of solving (Vapnik and Lerner [145])

(51) $$\underset{w,b}{\operatorname{minimize}} \tfrac{1}{2}\|w\|^2 \quad \text{subject to} \quad y_i(\langle w, x \rangle + b) \geq 1.$$

Note that $\|w\|^{-1} f(x_i)$ is the distance of the point $x_i$ to the hyperplane $H(w, b) := \{x | \langle w, x \rangle + b = 0\}$. The condition $y_i f(x_i) \geq 1$ implies that the margin of separation is at least $2\|w\|^{-1}$. The bound becomes exact if equality is attained for some $y_i = 1$ and $y_j = -1$. Consequently, minimizing $\|w\|$ subject to the constraints maximizes the margin of separation. Equation (51) is a quadratic program which can be solved efficiently (Fletcher [51]).

Mangasarian [95] devised a similar optimization scheme using $\|w\|_1$ instead of $\|w\|_2$ in the objective function of (51). The result is a *linear* program. In general, one can show (Smola et al. [124]) that minimizing the $\ell_p$ norm of $w$ leads to the maximizing of the margin of separation in the $\ell_q$ norm where $\frac{1}{p} + \frac{1}{q} = 1$. The $\ell_1$ norm leads to sparse approximation schemes (see also Chen et al. [29]), whereas the $\ell_2$ norm can be extended to Hilbert spaces and kernels.

To deal with nonseparable problems, that is, cases when (51) is infeasible, we need to relax the constraints of the optimization problem. Bennett and Mangasarian [17] and Cortes and Vapnik [34] impose a linear penalty on the violation of the large-margin constraints to obtain

(52) $$\underset{w,b,\xi}{\operatorname{minimize}} \tfrac{1}{2}\|w\|^2 + C\sum_{i=1}^{n} \xi_i$$

$$\text{subject to} \quad y_i(\langle w, x_i \rangle + b) \geq 1 - \xi_i \text{ and } \xi_i \geq 0, \forall i \in [n].$$

Equation (52) is a quadratic program which is always feasible (e.g., $w, b = 0$ and $\xi_i = 1$ satisfy the constraints). $C > 0$ is a regularization constant trading off the violation of the constraints vs. maximizing the overall margin.

Whenever the dimensionality of $\mathcal{X}$ exceeds $n$, direct optimization of (52) is computationally inefficient. This is particularly true if we map from $\mathcal{X}$

22T. HOFMANN, B. SCHÖLKOPF AND A. J. SMOLAinto an RKHS. To address these problems, one may solve the problem in dual space as follows. The Lagrange function of (52) is given by

$$\mathcal{L}(w,b,\xi,\alpha,\eta) = \tfrac{1}{2}\|w\|^2 + C\sum_{i=1}^{n}\xi_i$$
$$(53)$$
$$+ \sum_{i=1}^{n}\alpha_i(1 - \xi_i - y_i(\langle w, x_i\rangle + b)) - \sum_{i=1}^{n}\eta_i\xi_i,$$

where $\alpha_i, \eta_i \geq 0$ for all $i \in [n]$. To compute the dual of $\mathcal{L}$, we need to identify the first order conditions in $w, b$. They are given by

$$\partial_w \mathcal{L} = w - \sum_{i=1}^{n}\alpha_i y_i x_i = 0 \quad \text{and}$$
$$(54) \qquad \partial_b \mathcal{L} = -\sum_{i=1}^{n}\alpha_i y_i = 0 \quad \text{and}$$
$$\partial_{\xi_i}\mathcal{L} = C - \alpha_i + \eta_i = 0.$$

This translates into $w = \sum_{i=1}^{n}\alpha_i y_i x_i$, the linear constraint $\sum_{i=1}^{n}\alpha_i y_i = 0$, and the box-constraint $\alpha_i \in [0, C]$ arising from $\eta_i \geq 0$. Substituting (54) into $\mathcal{L}$ yields the Wolfe dual

$$\underset{\alpha}{\text{minimize}}\, \tfrac{1}{2}\alpha^\top Q\alpha - \alpha^\top 1 \quad \text{subject to} \quad \alpha^\top y = 0 \text{ and } \alpha_i \in [0, C],\, \forall i \in [n].$$
(55)

$Q \in \mathbb{R}^{n \times n}$ is the matrix of inner products $Q_{ij} := y_i y_j \langle x_i, x_j\rangle$. Clearly, this can be extended to feature maps and kernels easily via $K_{ij} := y_i y_j \langle \Phi(x_i), \Phi(x_j)\rangle = y_i y_j k(x_i, x_j)$. Note that $w$ lies in the span of the $x_i$. This is an instance of the representer theorem (Theorem 9). The KKT conditions (Boser et al. [23], Cortes and Vapnik [34], Karush [81] and Kuhn and Tucker [88]) require that at optimality $\alpha_i(y_i f(x_i) - 1) = 0$. This means that only those $x_i$ may appear in the expansion (54) for which $y_i f(x_i) \leq 1$, as otherwise $\alpha_i = 0$. The $x_i$ with $\alpha_i > 0$ are commonly referred to as support vectors.

Note that $\sum_{i=1}^{n}\xi_i$ is an upper bound on the empirical risk, as $y_i f(x_i) \leq 0$ implies $\xi_i \geq 1$ (see also Lemma 10). The number of misclassified points $x_i$ itself depends on the configuration of the data and the value of $C$. Ben-David et al. [15] show that finding even an approximate minimum classification error solution is difficult. That said, it is possible to modify (52) such that a desired target number of observations violates $y_i f(x_i) \geq \rho$ for some $\rho \in \mathbb{R}$ by making the threshold itself a variable of the optimization problem (Schölkopf et al. [120]). This leads to the following optimization problem ($\nu$-SV classification):

$$\underset{w,b,\xi}{\text{minimize}}\, \tfrac{1}{2}\|w\|^2 + \sum_{i=1}^{n}\xi_i - n\nu\rho$$

(56)
$$\text{subject to} \quad y_i(\langle w, x_i \rangle + b) \geq \rho - \xi_i \text{ and } \xi_i \geq 0.$$

The dual of (56) is essentially identical to (55) with the exception of an additional constraint:

(57)
$$\underset{\alpha}{\text{minimize}} \tfrac{1}{2}\alpha^\top Q \alpha \quad \text{subject to} \quad \alpha^\top y = 0 \text{ and } \alpha^\top 1 = n\nu \text{ and } \alpha_i \in [0,1].$$

One can show that for every $C$ there exists a $\nu$ such that the solution of (57) is a multiple of the solution of (55). Schölkopf et al. [120] prove that solving (57) for which $\rho > 0$ satisfies the following:

1. $\nu$ is an upper bound on the fraction of margin errors.
2. $\nu$ is a lower bound on the fraction of SVs.

Moreover, under mild conditions, with probability 1, asymptotically, $\nu$ equals both the fraction of SVs and the fraction of errors.

This statement implies that whenever the data are sufficiently well separable (i.e., $\rho > 0$), $\nu$-SV classification finds a solution with a fraction of at most $\nu$ margin errors. Also note that, for $\nu = 1$, all $\alpha_i = 1$, that is, $f$ becomes an affine copy of the Parzen windows classifier (5).

3.2. *Estimating the support of a density.* We now extend the notion of linear separation to that of estimating the support of a density (Schölkopf et al. [117] and Tax and Duin [134]). Denote by $X = \{x_1, \ldots, x_n\} \subseteq \mathcal{X}$ the sample drawn from $\mathrm{P}(x)$. Let $\mathcal{C}$ be a class of measurable subsets of $\mathcal{X}$ and let $\lambda$ be a real-valued function defined on $\mathcal{C}$. The *quantile function* (Einmal and Mason [47]) with respect to $(\mathrm{P}, \lambda, \mathcal{C})$ is defined as

(58) $\quad U(\mu) = \inf\{\lambda(C) | \mathrm{P}(C) \geq \mu, C \in \mathcal{C}\} \quad \text{where } \mu \in (0,1].$

We denote by $C_\lambda(\mu)$ and $C_\lambda^m(\mu)$ the (not necessarily unique) $C \in \mathcal{C}$ that attain the infimum (when it is achievable) on $\mathrm{P}(x)$ and on the empirical measure given by $X$ respectively. A common choice of $\lambda$ is the Lebesgue measure, in which case $C_\lambda(\mu)$ is the minimum volume set $C \in \mathcal{C}$ that contains at least a fraction $\mu$ of the probability mass.

Support estimation requires us to find some $C_\lambda^m(\mu)$ such that $|\mathrm{P}(C_\lambda^m(\mu)) - \mu|$ is small. This is where the complexity trade-off enters: On the one hand, we want to use a rich class $\mathcal{C}$ to capture all possible distributions, on the other hand, large classes lead to large deviations between $\mu$ and $\mathrm{P}(C_\lambda^m(\mu))$. Therefore, we have to consider classes of sets which are suitably restricted. This can be achieved using an SVM regularizer.

SV support estimation works by using SV support estimation related to previous work as follows: set $\lambda(C_w) = \|w\|^2$, where $C_w = \{x | f_w(x) \geq \rho\}$, $f_w(x) = \langle w, x \rangle$, and $(w, \rho)$ are respectively a weight vector and an offset.

Stated as a convex optimization problem, we want to separate the data from the origin with maximum margin via

$$\operatorname*{minimize}_{w,\xi,\rho} \tfrac{1}{2}\|w\|^2 + \sum_{i=1}^n \xi_i - n\nu\rho$$

(59)

$$\text{subject to} \quad \langle w, x_i \rangle \geq \rho - \xi_i \text{ and } \xi_i \geq 0.$$

Here, $\nu \in (0,1]$ plays the same role as in (56), controlling the number of observations $x_i$ for which $f(x_i) \leq \rho$. Since nonzero slack variables $\xi_i$ are penalized in the objective function, if $w$ and $\rho$ solve this problem, then the decision function $f(x)$ will attain or exceed $\rho$ for at least a fraction $1 - \nu$ of the $x_i$ contained in $X$, while the regularization term $\|w\|$ will still be small. The dual of (59) yield:

(60) $\quad \operatorname*{minimize}_{\alpha} \tfrac{1}{2}\alpha^\top K \alpha \quad \text{subject to} \quad \alpha^\top 1 = \nu n \text{ and } \alpha_i \in [0,1].$

To compare (60) to a Parzen windows estimator, assume that $k$ is such that it can be normalized as a density in input space, such as a Gaussian. Using $\nu = 1$ in (60), the constraints automatically imply $\alpha_i = 1$. Thus, $f$ reduces to a Parzen windows estimate of the underlying density. For $\nu < 1$, the equality constraint (60) still ensures that $f$ is a thresholded density, now depending only on a *subset* of $X$—those which are important for deciding whether $f(x) \leq \rho$.

3.3. *Regression estimation.* SV regression was first proposed in Vapnik [140] and Vapnik et al. [144] using the so-called $\epsilon$-insensitive loss function. It is a direct extension of the soft-margin idea to regression: instead of requiring that $yf(x)$ exceeds some margin value, we now require that the values $y - f(x)$ are bounded by a margin on both sides. That is, we impose the soft constraints

(61) $\quad\quad\quad y_i - f(x_i) \leq \epsilon_i - \xi_i \quad \text{and} \quad f(x_i) - y_i \leq \epsilon_i - \xi_i^*,$

where $\xi_i, \xi_i^* \geq 0$. If $|y_i - f(x_i)| \leq \epsilon$, no penalty occurs. The objective function is given by the sum of the slack variables $\xi_i, \xi_i^*$ penalized by some $C > 0$ and a measure for the slope of the function $f(x) = \langle w, x \rangle + b$, that is, $\tfrac{1}{2}\|w\|^2$.

Before computing the dual of this problem, let us consider a somewhat more general situation where we use a range of different convex penalties for the deviation between $y_i$ and $f(x_i)$. One may check that minimizing $\tfrac{1}{2}\|w\|^2 + C \sum_{i=1}^m \xi_i + \xi_i^*$ subject to (61) is equivalent to solving

(62) $\quad \operatorname*{minimize}_{w,b,\xi} \tfrac{1}{2}\|w\|^2 + \sum_{i=1}^n \psi(y_i - f(x_i)) \quad \text{where } \psi(\xi) = \max(0, |\xi| - \epsilon).$

Choosing different loss functions $\psi$ leads to a rather rich class of estimators:

- $\psi(\xi) = \frac{1}{2}\xi^2$ yields penalized least squares (LS) regression (Hoerl and Kennard [68], Morozov [102], Tikhonov [136] and Wahba [148]). The corresponding optimization problem can be minimized by solving a linear system.
- For $\psi(\xi) = |\xi|$, we obtain the penalized least absolute deviations (LAD) estimator (Bloomfield and Steiger [20]). That is, we obtain a quadratic program to estimate the conditional median.
- A combination of LS and LAD loss yields a penalized version of Huber's robust regression (Huber [71] and Smola and Schölkopf [126]). In this case we have $\psi(\xi) = \frac{1}{2\sigma}\xi^2$ for $|\xi| \le \sigma$ and $\psi(\xi) = |\xi| - \frac{\sigma}{2}$ for $|\xi| \ge \sigma$.
- Note that also quantile regression can be modified to work with kernels (Schölkopf et al. [120]) by using as loss function the "pinball" loss, that is, $\psi(\xi) = (1-\tau)\psi$ if $\psi < 0$ and $\psi(\xi) = \tau\psi$ if $\psi > 0$.

All the optimization problems arising from the above five cases are convex quadratic programs. Their dual resembles that of (61), namely,

(63a) $\quad \underset{\alpha,\alpha^*}{\text{minimize}} \frac{1}{2}(\alpha - \alpha^*)^\top K(\alpha - \alpha^*) + \epsilon^\top(\alpha + \alpha^*) - y^\top(\alpha - \alpha^*)$

(63b) $\quad \text{subject to} \quad (\alpha - \alpha^*)^\top 1 = 0 \text{ and } \alpha_i, \alpha_i^* \in [0, C].$

Here $K_{ij} = \langle x_i, x_j \rangle$ for linear models and $K_{ij} = k(x_i, x_j)$ if we map $x \to \Phi(x)$. The $\nu$-trick, as described in (56) (Schölkopf et al. [120]), can be extended to regression, allowing one to choose the margin of approximation automatically. In this case (63a) drops the terms in $\epsilon$. In its place, we add a linear constraint $(\alpha - \alpha^*)^\top 1 = \nu n$. Likewise, LAD is obtained from (63) by dropping the terms in $\epsilon$ without additional constraints. Robust regression leaves (63) unchanged, however, in the definition of $K$ we have an additional term of $\sigma^{-1}$ on the main diagonal. Further details can be found in Schölkopf and Smola [118]. For quantile regression we drop $\epsilon$ and we obtain different constants $C(1-\tau)$ and $C\tau$ for the constraints on $\alpha^*$ and $\alpha$. We will discuss uniform convergence properties of the empirical risk estimates with respect to various $\psi(\xi)$ in Section 3.6.

3.4. *Multicategory classification, ranking and ordinal regression.* Many estimation problems cannot be described by assuming that $\mathcal{Y} = \{\pm 1\}$. In this case it is advantageous to go beyond simple functions $f(x)$ depending on $x$ only. Instead, we can encode a larger degree of information by estimating a function $f(x,y)$ and subsequently obtaining a prediction via $\hat{y}(x) := \arg\max_{y \in \mathcal{Y}} f(x,y)$. In other words, we study problems where $y$ is obtained as the solution of an optimization problem over $f(x,y)$ and we wish to find $f$ such that $y$ matches $y_i$ as well as possible for relevant inputs $x$.

Note that the loss may be more than just a simple 0–1 loss. In the following we denote by $\Delta(y, y')$ the loss incurred by estimating $y'$ instead of $y$. Without loss of generality, we require that $\Delta(y, y) = 0$ and that $\Delta(y, y') \geq 0$ for all $y, y' \in \mathcal{Y}$. Key in our reasoning is the following:

LEMMA 10. *Let $f : \mathcal{X} \times \mathcal{Y} \to \mathbb{R}$ and assume that $\Delta(y, y') \geq 0$ with $\Delta(y, y) = 0$. Moreover, let $\xi \geq 0$ such that $f(x, y) - f(x, y') \geq \Delta(y, y') - \xi$ for all $y' \in \mathcal{Y}$. In this case $\xi \geq \Delta(y, \arg\max_{y' \in \mathcal{Y}} f(x, y'))$.*

The construction of the estimator was suggested by Taskar et al. [132] and Tsochantaridis et al. [137], and a special instance of the above lemma is given by Joachims [78]. While the bound appears quite innocuous, it allows us to describe a much richer class of estimation problems as a convex program.

To deal with the added complexity, we assume that $f$ is given by $f(x, y) = \langle \Phi(x, y), w \rangle$. Given the possibly nontrivial connection between $x$ and $y$, the use of $\Phi(x, y)$ cannot be avoided. Corresponding kernel functions are given by $k(x, y, x', y') = \langle \Phi(x, y), \Phi(x', y') \rangle$. We have the following optimization problem (Tsochantaridis et al. [137]):

$$\underset{w, \xi}{\text{minimize}} \tfrac{1}{2} \|w\|^2 + C \sum_{i=1}^{n} \xi_i \tag{64}$$

$$\text{subject to} \quad \langle w, \Phi(x_i, y_i) - \Phi(x_i, y) \rangle \geq \Delta(y_i, y) - \xi_i, \forall i \in [n], y \in \mathcal{Y}.$$

This is a convex optimization problem which can be solved efficiently if the constraints can be evaluated without high computational cost. One typically employs column-generation methods (Bennett et al. [16], Fletcher [51], Hettich and Kortanek [66] and Tsochantaridis et al. [137]) which identify one violated constraint at a time to find an approximate minimum of the optimization problem.

To describe the flexibility of the framework set out by (64) we give several examples below:

- Binary classification can be recovered by setting $\Phi(x, y) = y\Phi(x)$, in which case the constraint of (64) reduces to $2y_i \langle \Phi(x_i), w \rangle \geq 1 - \xi_i$. Ignoring constant offsets and a scaling factor of 2, this is exactly the standard SVM optimization problem.
- Multicategory classification problems (Allwein et al. [2], Collins [30] and Crammer and Singer [35]) can be encoded via $\mathcal{Y} = [N]$, where $N$ is the number of classes and $\Delta(y, y') = 1 - \delta_{y, y'}$. In other words, the loss is 1 whenever we predict the wrong class and 0 for correct classification. Corresponding kernels are typically chosen to be $\delta_{y, y'} k(x, x')$.
- We can deal with joint labeling problems by setting $\mathcal{Y} = \{\pm 1\}^n$. In other words, the error measure does not depend on a single observation but on an entire set of labels. Joachims [78] shows that the so-called $F_1$ score (van Rijsbergen [138]) used in document retrieval and the area under the ROC curve (Bamber [10]) fall into this category of problems. Moreover, Joachims [78] derives an $O(n^2)$ method for evaluating the inequality constraint over $\mathcal{Y}$.

- Multilabel estimation problems deal with the situation where we want to find the best subset of labels $\mathcal{Y} \subseteq 2^{[N]}$ which correspond to some observation $x$. Elisseeff and Weston [48] devise a ranking scheme where $f(x,i) > f(x,j)$ if label $i \in y$ and $j \notin y$. It is a special case of an approach described next.

Note that (64) is invariant under translations $\Phi(x,y) \leftarrow \Phi(x,y) + \Phi_0$ where $\Phi_0$ is constant, as $\Phi(x_i, y_i) - \Phi(x_i, y)$ remains unchanged. In practice, this means that transformations $k(x, y, x', y') \leftarrow k(x, y, x', y') + \langle \Phi_0, \Phi(x, y) \rangle + \langle \Phi_0, \Phi(x', y') \rangle + \|\Phi_0\|^2$ do not affect the outcome of the estimation process. Since $\Phi_0$ was arbitrary, we have the following lemma:

LEMMA 11. *Let $\mathcal{H}$ be an RKHS on $\mathcal{X} \times \mathcal{Y}$ with kernel $k$. Moreover, let $g \in \mathcal{H}$. Then the function $k(x, y, x', y') + f(x, y) + f(x', y') + \|g\|_{\mathcal{H}}^2$ is a kernel and it yields the same estimates as $k$.*

We need a slight extension to deal with general ranking problems. Denote by $\mathcal{Y} = \text{Graph}[N]$ the set of all directed graphs on $N$ vertices which do not contain loops of less than three nodes. Here an edge $(i, j) \in y$ indicates that $i$ is preferred to $j$ with respect to the observation $x$. It is the goal to find some function $f : \mathcal{X} \times [N] \to \mathbb{R}$ which imposes a total order on $[N]$ (for a given $x$) by virtue of the function values $f(x, i)$ such that the total order and $y$ are in good agreement.

More specifically, Crammer and Singer [36] and Dekel et al. [45] propose a decomposition algorithm $\mathcal{A}$ for the graphs $y$ such that the estimation error is given by the number of subgraphs of $y$ which are in disagreement with the total order imposed by $f$. As an example, multiclass classification can be viewed as a graph $y$ where the correct label $i$ is at the root of a directed graph and all incorrect labels are its children. Multilabel classification is then a bipartite graph where the correct labels only contain outgoing arcs and the incorrect labels only incoming ones.

This setting leads to a form similar to (64), except for the fact that we now have constraints over each subgraph $G \in \mathcal{A}(y)$. We solve

$$\underset{w,\xi}{\text{minimize}} \tfrac{1}{2}\|w\|^2 + C\sum_{i=1}^{n} |\mathcal{A}(y_i)|^{-1} \sum_{G \in \mathcal{A}(y_i)} \xi_{iG}$$

(65) $\quad\text{subject to}\quad \langle w, \Phi(x_i, u) - \Phi(x_i, v) \rangle \geq 1 - \xi_{iG}$ and $\xi_{iG} \geq 0$

$$\text{for all } (u, v) \in G \in \mathcal{A}(y_i).$$

That is, we test for all $(u, v) \in G$ whether the ranking imposed by $G \in y_i$ is satisfied.

Finally, ordinal regression problems which perform ranking not over labels $y$ but rather over observations $x$ were studied by Herbrich et al. [65] and Chapelle and Harchaoui [27] in the context of ordinal regression and conjoint analysis respectively. In ordinal regression $x$ is preferred to $x'$ if $f(x) > f(x')$ and, hence, one minimizes an optimization problem akin to (64), with constraint $\langle w, \Phi(x_i) - \Phi(x_j) \rangle \geq 1 - \xi_{ij}$. In conjoint analysis the same operation is carried out for $\Phi(x, u)$, where $u$ is the user under consideration. Similar models were also studied by Basilico and Hofmann [13]. Further models will be discussed in Section 4, in particular situations where $\mathcal{Y}$ is of exponential size. These models allow one to deal with sequences and more sophisticated structures.

3.5. *Applications of SVM algorithms.* When SVMs were first presented, they initially met with skepticism in the statistical community. Part of the reason was that, as described, SVMs construct their decision rules in potentially very high-dimensional feature spaces associated with kernels. Although there was a fair amount of theoretical work addressing this issue (see Section 3.6 below), it was probably to a larger extent the empirical success of SVMs that paved its way to become a standard method of the statistical toolbox. The first successes of SVMs on practical problems were in handwritten digit recognition, which was the main benchmark task considered in the Adaptive Systems Department at AT&T Bell Labs where SVMs were developed. Using methods to incorporate transformation invariances, SVMs were shown to beat the world record on the MNIST benchmark set, at the time the gold standard in the field (DeCoste and Schölkopf [44]). There has been a significant number of further computer vision applications of SVMs since then, including tasks such as object recognition and detection. Nevertheless, it is probably fair to say that two other fields have been more influential in spreading the use of SVMs: bioinformatics and natural language processing. Both of them have generated a spectrum of challenging high-dimensional problems on which SVMs excel, such as microarray processing tasks and text categorization. For references, see Joachims [77] and Schölkopf et al. [121].

Many successful applications have been implemented using SV classifiers; however, also the other variants of SVMs have led to very good results, including SV regression, SV novelty detection, SVMs for ranking and, more recently, problems with interdependent labels (McCallum et al. [96] and Tsochantaridis et al. [137]).

At present there exists a large number of readily available software packages for SVM optimization. For instance, SVMStruct, based on Tsochantaridis et al. [137] solves structured estimation problems. LibSVM is an open source solver which excels on binary problems. The Torch package contains a number of estimation methods, including SVM solvers. Several SVM implementations are also available via statistical packages, such as R.

3.6. *Margins and uniform convergence bounds.* While the algorithms were motivated by means of their practicality and the fact that 0–1 loss functions yield hard-to-control estimators, there exists a large body of work on statistical analysis. We refer to the works of Bartlett and Mendelson [12], Jordan et al. [80], Koltchinskii [86], Mendelson [98] and Vapnik [141] for details. In particular, the review of Bousquet et al. [24] provides an excellent summary of the current state of the art. Specifically for the structured case, recent work by Collins [30] and Taskar et al. [132] deals with explicit constructions to obtain better scaling behavior in terms of the number of class labels.

The general strategy of the analysis can be described by the following three steps: first, the discrete loss is upper bounded by some function, such as $\psi(yf(x))$, which can be efficiently minimized [e.g. the soft margin function $\max(0, 1 - yf(x))$ of the previous section satisfies this property]. Second, one proves that the empirical average of the $\psi$-loss is concentrated close to its expectation. This will be achieved by means of Rademacher averages. Third, one shows that under rather general conditions the minimization of the $\psi$-loss is consistent with the minimization of the expected risk. Finally, these bounds are combined to obtain rates of convergence which only depend on the Rademacher average and the approximation properties of the function class under consideration.

**4. Statistical models and RKHS.** As we have argued so far, the reproducing kernel Hilbert space approach offers many advantages in machine learning: (i) powerful and flexible models can be defined, (ii) many results and algorithms for linear models in Euclidean spaces can be generalized to RKHS, (iii) learning theory assures that effective learning in RKHS is possible, for instance, by means of regularization.

In this chapter we will show how kernel methods can be utilized in the context of *statistical models*. There are several reasons to pursue such an avenue. First of all, in conditional modeling, it is often insufficient to compute a prediction without assessing confidence and reliability. Second, when dealing with multiple or structured responses, it is important to model *dependencies between responses* in addition to the dependence on a set of covariates. Third, incomplete data, be it due to missing variables, incomplete training targets or a model structure involving latent variables, needs to be dealt with in a principled manner. All of these issues can be addressed by using the RKHS approach to define statistical models and by combining kernels with statistical approaches such as exponential models, generalized linear models and Markov networks.

### 4.1. *Exponential RKHS models.*

4.1.1. *Exponential models.* Exponential models or *exponential families* are among the most important class of parametric models studied in statistics. Given a canonical vector of statistics $\Phi$ and a $\sigma$-finite measure $\nu$ over the sample space $\mathscr{X}$, an exponential model can be defined via its probability density with respect to $\nu$ (cf. Barndorff-Nielsen [11]),

$$
(66) \quad p(x;\theta) = \exp[\langle \theta, \Phi(x) \rangle - g(\theta)]
$$

$$
\text{where } g(\theta) := \ln \int_{\mathcal{X}} e^{\langle \theta, \Phi(x) \rangle} \, d\nu(x).
$$

The $m$-dimensional vector $\theta \in \Theta$ with $\Theta := \{\theta \in \mathbb{R}^m : g(\theta) < \infty\}$ is also called the *canonical parameter* vector. In general, there are multiple exponential representations of the same model via canonical parameters that are affinely related to one another (Murray and Rice [103]). A representation with minimal $m$ is called a minimal representation, in which case $m$ is the *order* of the exponential model. One of the most important properties of exponential families is that they have sufficient statistics of fixed dimensionality, that is, the joint density for i.i.d. random variables $X_1, X_2, \ldots, X_n$ is also exponential, the corresponding canonical statistics simply being $\sum_{i=1}^{n} \Phi(X_i)$. It is well known that much of the structure of exponential models can be derived from the log partition function $g(\theta)$, in particular,

$$
(67) \quad \nabla_\theta g(\theta) = \mu(\theta) := \mathbf{E}_\theta[\Phi(X)], \qquad \partial_\theta^2 g(\theta) = \mathbf{V}_\theta[\Phi(X)],
$$

where $\mu$ is known as the mean-value map. Being a covariance matrix, the Hessian of $g$ is positive semi-definite and, consequently, $g$ is convex.

Maximum likelihood estimation (MLE) in exponential families leads to a particularly elegant form for the MLE equations: the expected and the observed canonical statistics agree at the MLE $\hat{\theta}$. This means, given an i.i.d. sample $\mathcal{S} = (x_i)_{i \in [n]}$,

$$
(68) \quad \mathbf{E}_{\hat{\theta}}[\Phi(X)] = \mu(\hat{\theta}) = \frac{1}{n} \sum_{i=1}^{n} \Phi(x_i) := r\mathbf{E}_{\mathcal{S}}[\Phi(X)].
$$

4.1.2. *Exponential RKHS models.* One can extend the parameteric exponential model in (66) by defining a statistical model via an RKHS $\mathcal{H}$ with generating kernel $k$. Linear function $\langle \theta, \Phi(\cdot) \rangle$ over $\mathcal{X}$ are replaced with functions $f \in \mathcal{H}$, which yields an exponential RKHS model

$$
(69) \quad p(x;f) = \exp[f(x) - g(f)],
$$

$$
f \in \mathcal{H} := \left\{ f : f(\cdot) = \sum_{x \in \mathcal{S}} \alpha_x k(\cdot, x), \mathcal{S} \subseteq \mathcal{X}, |\mathcal{S}| < \infty \right\}.
$$

A justification for using exponential RKHS families with rich canonical statistics as a generic way to define nonparametric models stems from the fact that if the chosen kernel $k$ is powerful enough, the associated exponential families become universal density estimators. This can be made precise using the concept of universal kernels (Steinwart [128], cf. Section 2).

PROPOSITION 12 (Dense densities). *Let $\mathcal{X}$ be a measurable set with a fixed $\sigma$-finite measure $\nu$ and denote by $\mathcal{P}$ a family of densities on $\mathcal{X}$ with respect to $\nu$ such that $p \in \mathcal{P}$ is uniformly bounded from above and continuous. Let $k : \mathcal{X} \times \mathcal{X} \to \mathbb{R}$ be a universal kernel for $\mathcal{H}$. Then the exponential RKHS family of densities generated by $k$ according to equation (69) is dense in $\mathcal{P}$ in the $L_\infty$ sense.*

4.1.3. *Conditional exponential models.* For the rest of the paper we will focus on the case of predictive or conditional modeling with a—potentially compound or structured—response variable $Y$ and predictor variables $X$. Taking up the concept of joint kernels introduced in the previous section, we will investigate conditional models that are defined by functions $f : \mathcal{X} \times \mathcal{Y} \to \mathbb{R}$ from some RKHS $\mathcal{H}$ over $\mathcal{X} \times \mathcal{Y}$ with kernel $k$ as follows:

$$p(y|x; f) = \exp[f(x,y) - g(x,f)]$$
(70)
$$\text{where } g(x, f) := \ln \int_{\mathcal{Y}} e^{f(x,y)} \, d\nu(y).$$

Notice that in the finite-dimensional case we have a feature map $\Phi : \mathcal{X} \times \mathcal{Y} \to \mathbb{R}^m$ from which parametric models are obtained via $\mathcal{H} := \{f : \exists w, f(x,y) = f(x,y;w) := \langle w, \Phi(x,y) \rangle\}$ and each $f$ can be identified with its parameter $w$. Let us discuss some concrete examples to illustrate the rather general model equation (70):

- Let $Y$ be univariate and define $\Phi(x, y) = y\Phi(x)$. Then simply $f(x, y; w) = \langle w, \Phi(x, y) \rangle = y\tilde{f}(x; w)$, with $\tilde{f}(x; w) := \langle w, \Phi(x) \rangle$ and the model equation in (70) reduces to

(71) $$p(y|x; w) = \exp[y\langle w, \Phi(x) \rangle - g(x, w)].$$

  This is a *generalized linear model* (GLM) (McCullagh and Nelder [97]) with a canonical link, that is, the canonical parameters depend linearly on the covariates $\Phi(x)$. For different response scales, we get several well-known models such as, for instance, logistic regression where $y \in \{-1, 1\}$.
- In the nonparameteric extension of generalized linear models following Green and Yandell [57] and O'Sullivan [105] the parametric assumption on the linear predictor $\tilde{f}(x; w) = \langle w, \Phi(x) \rangle$ in the GLMs is relaxed by requiring that $\tilde{f}$ comes from some sufficiently smooth class of functions, namely, an RKHS defined over $\mathcal{X}$. In combination with a parametric part, this can also be used to define semi-parametric models. Popular choices of kernels include the ANOVA kernel investigated by [149]. This is a special case of defining joint kernels from an existing kernel $k$ over inputs via $k((x,y),(x',y')) := yy'k(x,x')$.

- Joint kernels provide a powerful framework for prediction problems with structured outputs. An illuminating example is statistical natural language parsing with lexicalized probabilistic context free grammars (Magerman [94]). Here $x$ will be an English sentence and $y$ a parse tree for $x$, that is, a highly structured and complex output. The productions of the grammar are known, but the conditional probability $p(y|x)$ needs to be estimated based on training data of parsed/annotated sentences. In the simplest case, the extracted statistics $\Phi$ may encode the frequencies of the use of different productions in a sentence with a known parse tree. More sophisticated feature encodings are discussed in Taskar et al. [133] and Zettlemoyer and Collins [156]. The conditional modeling approach provide alternatives to state-of-the art approaches that estimate joint models $p(x,y)$ with maximum likelihood or maximum entropy and obtain predictive models by conditioning on $x$.

4.1.4. *Risk functions for model fitting.* There are different inference principles to determine the optimal function $f \in \mathcal{H}$ for the conditional exponential model in (70). One standard approach to parametric model fitting is to maximize the conditional log-likelihood—or equivalently—minimize a logarithmic loss, a strategy pursued in the Conditional Random Field (CRF) approach of Lafferty [90]. Here we consider the more general case of minimizing a functional that includes a monotone function of the Hilbert space norm $\|f\|_{\mathcal{H}}$ as a stabilizer (Wahba [148]). This reduces to penalized log-likelihood estimation in the finite-dimensional case,

$$C^{\mathrm{ll}}(f;\mathcal{S}) := -\frac{1}{n}\sum_{i=1}^{n}\ln p(y_i|x_i;f),$$

(72)

$$\hat{f}^{\mathrm{ll}}(\mathcal{S}) := \underset{f \in \mathcal{H}}{\arg\min}\, \frac{\lambda}{2}\|f\|_{\mathcal{H}}^2 + C^{\mathrm{ll}}(f;\mathcal{S}).$$

- For the parametric case, Lafferty et al. [90] have employed variants of improved iterative scaling (Darroch and Ratcliff [40] and Della Pietra [46]) to optimize equation (72), whereas Sha and Pereira [122] have investigated preconditioned conjugate gradient descent and limited memory quasi-Newton methods.
- In order to optimize equation (72) one usually needs to compute expectations of the canonical statistics $\mathbf{E}_f[\Phi(Y,x)]$ at sample points $x = x_i$, which requires the availability of efficient inference algorithms.

As we have seen in the case of classification and regression, likelihood-based criteria are by no means the only justifiable choice and large margin methods offer an interesting alternative. To that extend, we will present a general formulation of large margin methods for response variables over finite sample spaces that is based on the approach suggested by Altun et al. [6] and Taskar et al. [132]. Define

$$r(x, y; f) := f(x, y) - \max_{y' \neq y} f(x, y') = \min_{y' \neq y} \log \frac{p(y|x; f)}{p(y'|x; f)} \quad \text{and}$$
(73)
$$r(\mathcal{S}; f) := \min_{i=1}^{n} r(x_i, y_i; f).$$

Here $r(\mathcal{S}; f)$ generalizes the notion of separation margin used in SVMs. Since the log-odds ratio is sensitive to rescaling of $f$, that is, $r(x, y; \beta f) = \beta r(x, y; f)$, we need to constrain $\|f\|_{\mathcal{H}}$ to make the problem well defined. We thus replace $f$ by $\phi^{-1} f$ for some fixed dispersion parameter $\phi > 0$ and define the maximum margin problem $\hat{f}^{\mathrm{mm}}(\mathcal{S}) := \phi^{-1} \arg\max_{\|f\|_{\mathcal{H}}=1} r(\mathcal{S}; f/\phi)$. For the sake of the presentation, we will drop $\phi$ in the following. (We will not deal with the problem of how to estimate $\phi$ here; note, however, that one does need to know $\phi$ in order to make an optimal deterministic prediction.) Using the same line of arguments as was used in Section 3, the maximum margin problem can be re-formulated as a constrained optimization problem

(74) $$\hat{f}^{\mathrm{mm}}(\mathcal{S}) := \arg\min_{f \in \mathcal{H}} \tfrac{1}{2}\|f\|_{\mathcal{H}}^2 \quad \text{s.t.} \quad r(x_i, y_i; f) \geq 1, \forall i \in [n],$$

provided the latter is feasible, that is, if there exists $f \in \mathcal{H}$ such that $r(\mathcal{S}; f) > 0$. To make the connection to SVMs, consider the case of binary classification with $\Phi(x, y) = y\Phi(x)$, $f(x, y; w) = \langle w, y\Phi(x) \rangle$, where $r(x, y; f) = \langle w, y\Phi(x) \rangle - \langle w, -y\Phi(x) \rangle = 2y\langle w, \Phi(x) \rangle = 2\rho(x, y; w)$. The latter is twice the standard margin for binary classification in SVMs.

A soft margin version can be defined based on the Hinge loss as follows:

$$C^{\mathrm{hl}}(f; \mathcal{S}) := \frac{1}{n} \sum_{i=1}^{n} \min\{1 - r(x_i, y_i; f), 0\},$$
(75)
$$\hat{f}^{\mathrm{sm}}(\mathcal{S}) := \arg\min_{f \in \mathcal{H}} \frac{\lambda}{2}\|f\|_{\mathcal{H}}^2 + C^{\mathrm{hl}}(f, \mathcal{S}).$$

- An equivalent formulation using slack variables $\xi_i$ as discussed in Section 3 can be obtained by introducing soft-margin constraints $r(x_i, y_i; f) \geq 1 - \xi_i$, $\xi_i \geq 0$ and by defining $C^{\mathrm{hl}} = \frac{1}{n}\xi_i$. Each nonlinear constraint can be further expanded into $|\mathcal{Y}|$ linear constraints $f(x_i, y_i) - f(x_i, y) \geq 1 - \xi_i$ for all $y \neq y_i$.

- Prediction problems with structured outputs often involve task-specific loss function $\triangle : \mathcal{Y} \times \mathcal{Y} \to \mathbb{R}$ discussed in Section 3.4. As suggested in Taskar et al. [132] cost sensitive large margin methods can be obtained by defining re-scaled margin constraints $f(x_i, y_i) - f(x_i, y) \geq \triangle(y_i, y) - \xi_i$.
- Another sensible option in the parametric case is to minimize an exponential risk function of the following type:

$$\hat{f}^{\text{exp}}(\mathcal{S}) := \arg\min_w \frac{1}{n} \sum_{i=1}^n \sum_{y \neq y_i} \exp[f(x_i, y_i; w) - f(x_i, y; w)]. \tag{76}$$

This is related to the exponential loss used in the AdaBoost method of Freund and Schapire [53]. Since we are mainly interested in kernel-based methods here, we refrain from further elaborating on this connection.

4.1.5. *Generalized representer theorem and dual soft-margin formulation.* It is crucial to understand how the representer theorem applies in the setting of arbitrary discrete output spaces, since a finite representation for the optimal $\hat{f} \in \{\hat{f}^{\text{ll}}, \hat{f}^{\text{sm}}\}$ is the basis for constructive model fitting. Notice that the regularized log-loss, as well as the soft margin functional introduced above, depends not only on the values of $f$ on the sample $\mathcal{S}$, but rather on the evaluation of $f$ on the augmented sample $\tilde{\mathcal{S}} := \{(x_i, y) : i \in [n], y \in \mathcal{Y}\}$. This is the case, because for each $x_i$, output values $y \neq y_i$ not observed with $x_i$ show up in the log-partition function $g(x_i, f)$ in (70), as well as in the log-odds ratios in (73). This adds an additional complication compared to binary classification.

COROLLARY 13. *Denote by $\mathcal{H}$ an RKHS on $\mathcal{X} \times \mathcal{Y}$ with kernel $k$ and let $\mathcal{S} = ((x_i, y_i))_{i \in [n]}$. Furthermore, let $C(f; \mathcal{S})$ be a functional depending on $f$ only via its values on the augmented sample $\tilde{\mathcal{S}}$. Let $\Omega$ be a strictly monotonically increasing function. Then the solution of the optimization problem $\hat{f}(\mathcal{S}) := \arg\min_{f \in \mathcal{H}} C(f; \tilde{\mathcal{S}}) + \Omega(\|f\|_\mathcal{H})$ can be written as*

$$\hat{f}(\cdot) = \sum_{i=1}^n \sum_{y \in \mathcal{Y}} \beta_{iy} k(\cdot, (x_i, y)). \tag{77}$$

This follows directly from Theorem 9.

Let us focus on the soft margin maximizer $\hat{f}^{\text{sm}}$. Instead of solving (75) directly, we first derive the dual program, following essentially the derivation in Section 3.

PROPOSITION 14 (Tsochantaridis et al. [137]). *The minimizer $\hat{f}^{\text{sm}}(\mathcal{S})$ can be written as in Corollary 13, where the expansion coefficients can be* *computed from the solution of the following convex quadratic program:*

$$\text{(78a)} \quad \alpha^* = \arg\min_{\alpha}\left\{\tfrac{1}{2}\sum_{i,j=1}^{n}\sum_{y\neq y_i}\sum_{y'\neq y_j}\alpha_{iy}\alpha_{jy'}K_{iy,jy'} - \sum_{i=1}^{n}\sum_{y\neq y_i}\alpha_{iy}\right\}$$

$$\text{(78b)} \quad s.t. \quad \lambda n\sum_{y\neq y_i}\alpha_{iy} \leq 1,\ \forall i\in[n]; \alpha_{iy}\geq 0, \forall i\in[n], y\in\mathcal{Y},$$

*where* $K_{iy,jy'} := k((x_i,y_i),(x_j,y_j)) + k((x_i,y),(x_j,y')) - k((x_i,y_i),(x_j,y')) - k((x_i,y),(x_j,y_j)).$

- The multiclass SVM formulation of [35] can be recovered as a special case for kernels that are diagonal with respect to the outputs, that is, $k((x,y),(x',y')) = \delta_{y,y'}k(x,x')$. Notice that in this case the quadratic part in equation (78a) simplifies to

$$\sum_{i,j}k(x_i,x_j)\sum_{y}\alpha_{iy}\alpha_{jy}[1 + \delta_{y_i,y}\delta_{y_j,y} - \delta_{y_i,y} - \delta_{y_j,y}].$$

- The pairs $(x_i,y)$ for which $\alpha_{iy} > 0$ are the *support pairs*, generalizing the notion of support vectors. As in binary SVMs, their number can be much smaller than the total number of constraints. Notice also that in the final expansion contributions $k(\cdot,(x_i,y_i))$ will get nonnegative weights, whereas $k(\cdot,(x_i,y))$ for $y \neq y_i$ will get nonpositive weights. Overall one gets a balance equation $\beta_{iy_i} - \sum_{y\neq y_i}\beta iy = 0$ for every data point.

4.1.6. *Sparse approximation.* Proposition 14 shows that sparseness in the representation of $\hat{f}^{\text{sm}}$ is linked to the fact that only few $\alpha_{iy}$ in the solution to the dual problem in equation (78) are nonzero. Note that each of these Lagrange multipliers is linked to the corresponding soft margin constraints $f(x_i,y_i) - f(x_i,y) \geq 1 - \xi_i$. Hence, sparseness is achieved, if only few constraints are active at the optimal solution. While this may or may not be the case for a given sample, one can still exploit this observation to define a nested sequence of relaxations, where margin constraint are incrementally added. This corresponds to a constraint selection algorithm (Bertsimas and Tsitsiklis [19]) for the primal or, equivalently, a variable selection or column generation method for the dual program and has been investigated in Tsochantaridis et al. [137]. Solving a sequence of increasingly tighter relaxations to a mathematical problem is also known as an *outer approximation*. In particular, one may iterate through the training examples according to some (fair) visitation schedule and greedily select constraints that are most violated at the current solution $f$, that is, for the $i$th instance one computes

$$\text{(79)} \quad \hat{y}_i = \arg\max_{y\neq y_i} f(x_i,y) = \arg\max_{y\neq y_i} p(y|x_i;f),$$

and then strengthens the current relaxation by including $\alpha_{i\hat{y}_i}$ in the optimization of the dual if $f(x_i, y_i) - f(x_i, \hat{y}_i) < 1 - \xi_i - \epsilon$. Here $\epsilon > 0$ is a pre-defined tolerance parameter. It is important to understand how many strengthening steps are necessary to achieve a reasonable close approximation to the original problem. The following theorem provides an answer:

THEOREM 15 (Tsochantaridis et al. [137]). *Let $\bar{R} = \max_{i,y} K_{iy,iy}$ and choose $\epsilon > 0$. A sequential strengthening procedure, which optimizes equation (75) by greedily selecting $\epsilon$-violated constraints, will find an approximate solution where all constraints are fulfilled within a precision of $\epsilon$, that is, $r(x_i, y_i; f) \geq 1 - \xi_i - \epsilon$ after at most $\frac{2n}{\epsilon} \cdot \max\{1, \frac{4\bar{R}^2}{\lambda n^2 \epsilon}\}$ steps.*

COROLLARY 16. *Denote by $(\hat{f}, \hat{\xi})$ the optimal solution of a relaxation of the problem in Proposition 14, minimizing $\mathcal{R}(f, \xi, \mathcal{S})$ while violating no constraint by more than $\epsilon$ (cf. Theorem 15). Then*

$$\mathcal{R}(\hat{f}, \hat{\xi}, \mathcal{S}) \leq \mathcal{R}(\hat{f}^{\mathrm{sm}}, \xi^*, \mathcal{S}) \leq \mathcal{R}(\hat{f}, \hat{\xi}, \mathcal{S}) + \epsilon,$$

*where $(\hat{f}^{\mathrm{sm}}, \xi^*)$ is the optimal solution of the original problem.*

- Combined with an efficient QP solver, the above theorem guarantees a runtime polynomial in $n$, $\epsilon^{-1}$, $\bar{R}$ and $\lambda^{-1}$. This holds irrespective of special properties of the data set utilized, the only exception being the dependency on the sample points $x_i$ is through the radius $\bar{R}$.
- The remaining key problem is how to compute equation (79) efficiently. The answer depends on the specific form of the joint kernel $k$ and/or the feature map $\Phi$. In many cases, efficient dynamic programming techniques exists, whereas in other cases one has to resort to approximations or use other methods to identify a set of candidate distractors $\mathcal{Y}_i \subset \mathcal{Y}$ for a training pair $(x_i, y_i)$ (Collins [30]). Sometimes one may also have search heuristics available that may not find the solution to Equation (79), but that find (other) $\epsilon$-violating constraints with a reasonable computational effort.

4.1.7. *Generalized Gaussian processes classification.* The model equation (70) and the minimization of the regularized log-loss can be interpreted as a generalization of Gaussian process classification (Altun et al. [4] and Rasmussen and Williams [111]) by assuming that $(f(x, \cdot))_{x \in \mathcal{X}}$ is a vector-valued zero mean Gaussian process; note that the covariance function $C$ is defined over pairs $\mathcal{X} \times \mathcal{Y}$. For a given sample $\mathcal{S}$, define a multi-index vector $F(\mathcal{S}) := (f(x_i, y))_{i,y}$ as the restriction of the stochastic process $f$ to the augmented sample $\tilde{\mathcal{S}}$. Denote the kernel matrix by $K = (K_{iy,jy'})$, where $K_{iy,jy'} := C((x_i, y), (x_j, y'))$ with indices $i, j \in [n]$ and $y, y' \in \mathcal{Y}$, so that, in summary, $F(\mathcal{S}) \sim \mathcal{N}(0,K)$. This induces a predictive model via Bayesian model integration according to

$$(80) \qquad p(y|x;\mathcal{S}) = \int p(y|F(x,\cdot))p(F|\mathcal{S})\,dF,$$

where $x$ is a test point that has been included in the sample (transductive setting). For an i.i.d. sample, the log-posterior for $F$ can be written as

$$(81) \qquad \ln p(F|\mathcal{S}) = -\tfrac{1}{2}F^T\mathbf{K}^{-1}F + \sum_{i=1}^n [f(x_i,y_i) - g(x_i,F)] + const.$$

Invoking the representer theorem for $\hat{F}(\mathcal{S}) := \arg\max_F \ln p(F|\mathcal{S})$, we know that

$$(82) \qquad \hat{F}(\mathcal{S})_{iy} = \sum_{j=1}^n \sum_{y'\in\mathcal{Y}} \alpha_{iy} K_{iy,jy'},$$

which we plug into equation (81) to arrive at

$$(83) \qquad \min_\alpha \alpha^T\mathbf{K}\alpha - \sum_{i=1}^n \left( \alpha^T\mathbf{K}e_{iy_i} + \log \sum_{y\in\mathcal{Y}} \exp[\alpha^T\mathbf{K}e_{iy}] \right),$$

where $e_{iy}$ denotes the respective unit vector. Notice that for $f(\cdot) = \sum_{i,y} \alpha_{iy} k(\cdot, (x_i,y))$ the first term is equivalent to the squared RKHS norm of $f \in \mathcal{H}$ since $\langle f,f \rangle_\mathcal{H} = \sum_{i,j} \sum_{y,y'} \alpha_{iy}\alpha_{jy'} \langle k(\cdot,(x_i,y)), k(\cdot,(x_j,y')) \rangle$. The latter inner product reduces to $k((x_i,y),(x_j,y'))$ due to the reproducing property. Again, the key issue in solving (83) is how to achieve spareness in the expansion for $\hat{F}$.

4.2. *Markov networks and kernels.* In Section 4.1 no assumptions about the specific structure of the joint kernel defining the model in equation (70) has been made. In the following, we will focus on a more specific setting with multiple outputs, where dependencies are modeled by a conditional independence graph. This approach is motivated by the fact that independently predicting individual responses based on marginal response models will often be suboptimal and explicitly modeling these interactions can be of crucial importance.

4.2.1. *Markov networks and factorization theorem.* Denote predictor variables by $X$, response variables by $Y$ and define $Z := (X,Y)$ with associated sample space $\mathcal{Z}$. We use Markov networks as the modeling formalism for representing dependencies between covariates and response variables, as well as interdependencies among response variables.

DEFINITION 17. A *conditional independence graph* (or Markov network) is an undirected graph $\mathcal{G} = (Z, E)$ such that for any pair of variables $(Z_i, Z_j) \notin E$ if and only if $Z_i \perp\!\!\!\perp Z_j | Z - \{Z_i, Z_j\}$.

The above definition is based on the pairwise Markov property, but by virtue of the separation theorem (see, e.g., Whittaker [154]) this implies the global Markov property for distributions with full support. The global Markov property says that for disjoint subsets $U, V, W \subseteq Z$ where $W$ separates $U$ from $V$ in $\mathcal{G}$ one has that $U \perp\!\!\!\perp V | W$. Even more important in the context of this paper is the factorization result due to Hammersley and Clifford [61].

THEOREM 18. *Given a random vector $Z$ with conditional independence graph $\mathcal{G}$, any density function for $Z$ with full support factorizes over $\mathscr{C}(\mathcal{G})$, the set of maximal cliques of $\mathcal{G}$ as follows:*

$$p(z) = \exp\left[\sum_{c \in \mathscr{C}(G)} f_c(z_c)\right], \tag{84}$$

*where $f_c$ are clique compatibility functions dependent on $z$ only via the restriction on clique configurations $z_c$.*

The significance of this result is that in order to specify a distribution for $Z$, one only needs to specify or estimate the simpler functions $f_c$.

4.2.2. *Kernel decomposition over Markov networks.* It is of interest to analyze the structure of kernels $k$ that generate Hilbert spaces $\mathcal{H}$ of functions that are consistent with a graph.

DEFINITION 19. A function $f : \mathcal{Z} \to \mathbb{R}$ is *compatible* with a conditional independence graph $\mathcal{G}$, if $f$ decomposes additively as $f(z) = \sum_{c \in \mathscr{C}(\mathcal{G})} f_c(z_c)$ with suitably chosen functions $f_c$. A Hilbert space $\mathcal{H}$ over $\mathcal{Z}$ is compatible with $\mathcal{G}$, if every function $f \in \mathcal{H}$ is compatible with $\mathcal{G}$. Such $f$ and $\mathcal{H}$ are also called $\mathcal{G}$-compatible.

PROPOSITION 20. *Let $\mathcal{H}$ with kernel $k$ be a $\mathcal{G}$-compatible RKHS. Then there are functions $k_{cd} : \mathcal{Z}_c \times \mathcal{Z}_d \to \mathbb{R}$ such that the kernel decomposes as*

$$k(u, z) = \sum_{c,d \in \mathscr{C}} k_{cd}(u_c, z_d).$$

LEMMA 21. *Let $\mathcal{X}$ be a set of $n$-tupels and $f_i, g_i : \mathcal{X} \times \mathcal{X} \to \mathbb{R}$ for $i \in [n]$ functions such that $f_i(x, y) = f_i(x_i, y)$ and $g_i(x, y) = g_i(x, y_i)$. If $\sum_i f_i(x, y) = \sum_j g_j(x, y_j)$ for all $x, y$, then there exist functions $h_{ij}$ such that $\sum_i f_i(x, y) = \sum_{i,j} h_{ij}(x_i, y_j)$.*

- Proposition 20 is useful for the design of kernels, since it states that only kernels allowing an additive decomposition into local functions $k_{cd}$ are compatible with a given Markov network $\mathcal{G}$. Lafferty et al. [89] have pursued a similar approach by considering kernels for RKHS with functions defined over $\mathcal{Z}_{\mathscr{C}} := \{(c, z_c) : c \in c, z_c \in \mathcal{Z}_c\}$. In the latter case one can even deal with cases where the conditional dependency graph is (potentially) different for every instance.
- An illuminating example of how to design kernels via the decomposition in Proposition 20 is the case of *conditional Markov chains*, for which models based on joint kernels have been proposed in Altun et al. [6], Collins [30], Lafferty et al. [90] and Taskar et al. [132]. Given an input sequences $X = (X_t)_{t \in [T]}$, the goal is to predict a sequence of labels or class variables $Y = (Y_t)_{t \in [T]}, Y_t \in \Sigma$. Dependencies between class variables are modeled in terms of a Markov chain, whereas outputs $Y_t$ are assumed to depend (directly) on an observation window $(X_{t-r}, \ldots, X_t, \ldots, X_{t+r})$. Notice that this goes beyond the standard hidden Markov model structure by allowing for overlapping features ($r \geq 1$). For simplicity, we focus on a window size of $r = 1$, in which case the clique set is given by $\mathscr{C} := \{c_t := (x_t, y_t, y_{t+1}), c'_t := (x_{t+1}, y_t, y_{t+1}) : t \in [T-1]\}$. We assume an input kernel $k$ is given and introduce indicator vectors (or dummy variates) $I(Y_{\{t,t+1\}}) := (I_{\omega,\omega'}(Y_{\{t,t+1\}}))_{\omega,\omega' \in \Sigma}$. Now we can define the local kernel functions as

$$k_{cd}(z_c, z'_d) := \langle I(y_{\{s,s+1\}}) I(y'_{\{t,t+1\}}) \rangle$$

(85)
$$\times \begin{cases} k(x_s, x_t), & \text{if } c = c_s \text{ and } d = c_t, \\ k(x_{s+1}, x_{t+1}), & \text{if } c = c'_s \text{ and } d = c'_t. \end{cases}$$

Notice that the inner product between indicator vectors is zero, unless the variable pairs are in the same configuration.

Conditional Markov chain models have found widespread applications in natural language processing (e.g., for part of speech tagging and shallow parsing, cf. Sha and Pereira [122]), in information retrieval (e.g., for information extraction, cf. McCallum et al. [96]) or in computational biology (e.g., for gene prediction, cf. Culotta et al. [39]).

4.2.3. *Clique-based sparse approximation.* Proposition 20 immediately leads to an alternative version of the representer theorem as observed by Lafferty et al. [89] and Altum et al. [4].

COROLLARY 22. *If $\mathcal{H}$ is $\mathcal{G}$-compatible then in the same setting as in Corollary 13, the optimizer $\hat{f}$ can be written as*

(86) $$\hat{f}(u) = \sum_{i=1}^{n} \sum_{c \in \mathscr{C}} \sum_{y_c \in \mathcal{Y}_c} \beta^i_{c,y_c} \sum_{d \in \mathscr{C}} k_{cd}((x_{ic}, y_c), u_d),$$

here $x_{ic}$ are the variables of $x_i$ belonging to clique $c$ and $\mathcal{Y}_c$ is the subspace of $\mathcal{Z}_c$ that contains response variables.

- Notice that the number of parameters in the representation equation (86) scales with $n \cdot \sum_{c \in \mathscr{C}} |\mathcal{Y}_c|$ as opposed to $n \cdot |\mathcal{Y}|$ in equation (77). For cliques with reasonably small state spaces, this will be a significantly more compact representation. Notice also that the evaluation of functions $k_{cd}$ will typically be more efficient than evaluating $k$.
- In spite of this improvement, the number of terms in the expansion in equation (86) may in practice still be too large. In this case, one can pursue a reduced set approach, which selects a subset of variables to be included in a sparsified expansion. This has been proposed in Taskar et al. [132] for the soft margin maximization problem, as well as in Altun et al. [5] and Lafferty et al. [89] for conditional random fields and Gaussian processes. For instance, in Lafferty et al. [89] parameters $\beta^i_{cy_c}$ that maximize the functional gradient of the regularized log-loss are greedily included in the reduced set. In Taskar et al. [132] a similar selection criterion is utilized with respect to margin violations, leading to an SMO-like optimization algorithm (Platt [107]).

4.2.4. *Probabilistic inference.* In dealing with structured or interdependent response variables, computing marginal probabilities of interest or computing the most probable response [cf. equation (79)] may be nontrivial. However, for dependency graphs with small tree width, efficient inference algorithms exist, such as the junction tree algorithm (Dawid [43] and Jensen et al. [76]) and variants thereof. Notice that in the case of the conditional or hidden Markov chain, the junction tree algorithm is equivalent to the well-known forward–backward algorithm (Baum [14]). Recently, a number of approximate inference algorithms have been developed to deal with dependency graphs for which exact inference is not tractable (see, e.g., Wainwright and Jordan [150]).

**5. Kernel methods for unsupervised learning.** This section discusses various methods of data analysis by modeling the distribution of data in feature space. To that extent, we study the behavior of $\Phi(x)$ by means of rather simple linear methods, which have implications for nonlinear methods on the original data space $\mathcal{X}$. In particular, we will discuss the extension of PCA to Hilbert spaces, which allows for image denoising, clustering, and nonlinear dimensionality reduction, the study of covariance operators for the measure of independence, the study of mean operators for the design of two-sample tests, and the modeling of complex dependencies between sets of random variables via kernel dependency estimation and canonical correlation analysis.

5.1. *Kernel principal component analysis.* Principal component analysis (PCA) is a powerful technique for extracting structure from possibly high-dimensional data sets. It is readily performed by solving an eigenvalue problem, or by using iterative algorithms which estimate principal components.

PCA is an orthogonal transformation of the coordinate system in which we describe our data. The new coordinate system is obtained by projection onto the so-called principal axes of the data. A small number of principal components is often sufficient to account for most of the structure in the data.

The basic idea is strikingly simple: denote by $X = \{x_1, \ldots, x_n\}$ an $n$-sample drawn from P($x$). Then the covariance operator $C$ is given by $C = \mathbf{E}[(x - \mathbf{E}[x])(x - \mathbf{E}[x])^\top]$. PCA aims at estimating leading eigenvectors of $C$ via the empirical estimate $C_{\text{emp}} = \mathbf{E}_{\text{emp}}[(x - \mathbf{E}_{\text{emp}}[x])(x - \mathbf{E}_{\text{emp}}[x])^\top]$. If $\mathcal{X}$ is $d$-dimensional, then the eigenvectors can be computed in $O(d^3)$ time (Press et al. [110]).

The problem can also be posed in feature space (Schölkopf et al. [119]) by replacing $x$ with $\Phi(x)$. In this case, however, it is impossible to compute the eigenvectors directly. Yet, note that the image of $C_{\text{emp}}$ lies in the span of $\{\Phi(x_1), \ldots, \Phi(x_n)\}$. Hence, it is sufficient to diagonalize $C_{\text{emp}}$ in that subspace. In other words, we replace the *outer* product $C_{\text{emp}}$ by an *inner* product matrix, leaving the eigenvalues unchanged, which can be computed efficiently. Using $w = \sum_{i=1}^n \alpha_i \Phi(x_i)$, it follows that $\alpha$ needs to satisfy $PKP\alpha = \lambda\alpha$, where $P$ is the projection operator with $P_{ij} = \delta_{ij} - n^{-2}$ and $K$ is the kernel matrix on $X$.

Note that the problem can also be recovered as one of maximizing some Contrast$[f, X]$ subject to $f \in \mathcal{F}$. This means that the projections onto the leading eigenvectors correspond to the most reliable features. This optimization problem also allows us to unify various feature extraction methods as follows:

- For Contrast$[f, X] = \text{Var}_{\text{emp}}[f, X]$ and $\mathcal{F} = \{\langle w, x \rangle$ subject to $\|w\| \leq 1\}$, we recover PCA.
- Changing $\mathcal{F}$ to $\mathcal{F} = \{\langle w, \Phi(x) \rangle$ subject to $\|w\| \leq 1\}$, we recover kernel PCA.
- For Contrast$[f, X] = \text{Curtosis}[f, X]$ and $\mathcal{F} = \{\langle w, x \rangle$ subject to $\|w\| \leq 1\}$, we have Projection Pursuit (Friedman and Tukey [55] and Huber [72]). Other contrasts lead to further variants, that is, the Epanechikov kernel, entropic contrasts, and so on (Cook et al. [32], Friedman [54] and Jones and Sibson [79]).
- If $\mathcal{F}$ is a convex combination of basis functions and the contrast function is convex in $w$, one obtains computationally efficient algorithms, as the solution of the optimization problem can be found at one of the vertices (Rockafellar [114] and Schölkopf and Smola [118]).

Subsequent projections are obtained, for example, by seeking directions orthogonal to $f$ or other computationally attractive variants thereof.

Kernel PCA has been applied to numerous problems, from preprocessing and invariant feature extraction (Mika et al. [100]) to image denoising and super-resolution (Kim et al. [84]). The basic idea in the latter case is to obtain a set of principal directions in feature space $w_1, \ldots, w_l$, obtained from noise-free data, and to project the image $\Phi(x)$ of a noisy observation $x$ onto the space spanned by $w_1, \ldots, w_l$. This yields a "denoised" solution $\tilde{\Phi}(x)$ in feature space. Finally, to obtain the pre-image of this denoised solution, one minimizes $\|\Phi(x') - \tilde{\Phi}(x)\|$. The fact that projections onto the leading principal components turn out to be good starting points for pre-image iterations is further exploited in kernel dependency estimation (Section 5.3). Kernel PCA can be shown to contain several popular dimensionality reduction algorithms as special cases, including LLE, Laplacian Eigenmaps and (approximately) Isomap (Ham et al. [60]).

5.2. *Canonical correlation and measures of independence.* Given two samples $X, Y$, canonical correlation analysis (Hotelling [70]) aims at finding directions of projection $u, v$ such that the correlation coefficient between $X$ and $Y$ is maximized. That is, $(u, v)$ are given by

$$\underset{u,v}{\arg\max} \operatorname{Var}_{\text{emp}}[\langle u, x \rangle]^{-1} \operatorname{Var}_{\text{emp}}[\langle v, y \rangle]^{-1}$$
(87)
$$\times \mathbf{E}_{\text{emp}}[\langle u, x - \mathbf{E}_{\text{emp}}[x] \rangle \langle v, y - \mathbf{E}_{\text{emp}}[y] \rangle].$$

This problem can be solved by finding the eigensystem of $C_x^{-1/2} C_{xy} C_y^{-1/2}$, where $C_x, C_y$ are the covariance matrices of $X$ and $Y$ and $C_{xy}$ is the covariance matrix between $X$ and $Y$, respectively. Multivariate extensions are discussed in Kettenring [83].

CCA can be extended to kernels by means of replacing linear projections $\langle u, x \rangle$ by projections in feature space $\langle u, \Phi(x) \rangle$. More specifically, Bach and Jordan [8] used the so-derived contrast to obtain a measure of independence and applied it to Independent Component Analysis with great success. However, the formulation requires an additional regularization term to prevent the resulting optimization problem from becoming distribution independent.

Rényi [113] showed that independence between random variables is equivalent to the condition of vanishing covariance $\operatorname{Cov}[f(x), g(y)] = 0$ for all $C^1$ functions $f, g$ bounded by $L_\infty$ norm 1 on $\mathcal{X}$ and $\mathcal{Y}$. In Bach and Jordan [8], Das and Sen [41], Dauxois and Nkiet [42] and Gretton et al. [58, 59] a constrained empirical estimate of the above criterion is used. That is, one studies

$$\Lambda(X, Y, \mathcal{F}, \mathcal{G}) := \sup_{f,g} \operatorname{Cov}_{\text{emp}}[f(x), g(y)]$$
(88)
$$\text{subject to} \quad f \in \mathcal{F} \text{ and } g \in \mathcal{G}.$$

This statistic is often extended to use the entire series $\Lambda_1, \ldots, \Lambda_d$ of maximal correlations where each of the function pairs $(f_i, g_i)$ are orthogonal to the previous set of terms. More specifically Douxois and Nkiet [42] restrict $\mathcal{F}, \mathcal{G}$ to finite-dimensional linear function classes subject to their $L_2$ norm bounded by 1, Bach and Jordan [8] use functions in the RKHS for which some sum of the $\ell_2^n$ and the RKHS norm on the sample is bounded.

Gretton et al. [58] use functions with bounded RKHS norm only, which provides necessary and sufficient criteria if kernels are universal. That is, $\Lambda(X, Y, \mathcal{F}, \mathcal{G}) = 0$ if and only if $x$ and $y$ are independent. Moreover, $\operatorname{tr} PK_x PK_y P$ has the same theoretical properties and it can be computed much more easily in linear time, as it allows for incomplete Cholesky factorizations. Here $K_x$ and $K_y$ are the kernel matrices on $X$ and $Y$ respectively.

The above criteria can be used to derive algorithms for Independent Component Analysis (Bach and Jordan [8] and Gretton et al. [58]). While these algorithms come at a considerable computational cost, they offer very good performance. For faster algorithms, consider the work of Cardoso [26], Hyvärinen [73] and Lee et al. [91]. Also, the work of Chen and Bickel [28] and Yang and Amari [155] is of interest in this context.

Note that a similar approach can be used to develop two-sample tests based on kernel methods. The basic idea is that for universal kernels the map between distributions and points on the marginal polytope $\mu : p \to \mathbf{E}_{x \sim p}[\phi(x)]$ is bijective and, consequently, it imposes a norm on distributions. This builds on the ideas of [52]. The corresponding distance $d(p, q) := \|\mu[p] - \mu[q]\|$ leads to a $U$-statistic which allows one to compute empirical estimates of distances between distributions efficiently [22].

5.3. *Kernel dependency estimation.* A large part of the previous discussion revolved around estimating dependencies between samples $\mathcal{X}$ and $\mathcal{Y}$ for rather structured spaces $\mathcal{Y}$, in particular, (64). In general, however, such dependencies can be hard to compute. Weston et al. [153] proposed an algorithm which allows one to extend standard regularized LS regression models, as described in Section 3.3, to cases where $\mathcal{Y}$ has complex structure.

It works by recasting the estimation problem as a *linear* estimation problem for the map $f : \Phi(x) \to \Phi(y)$ and then as a nonlinear pre-image estimation problem for finding $\hat{y} := \operatorname{argmin}_y \|f(x) - \Phi(y)\|$ as the point in $\mathcal{Y}$ closest to $f(x)$.

This problem can be solved directly (Cortes et al. [33]) without the need for subspace projections. The authors apply it to the analysis of sequence data.

**6. Conclusion.** We have summarized some of the advances in the field of machine learning with positive definite kernels. Due to lack of space, this article is by no means comprehensive, in particular, we were not able to cover statistical learning theory, which is often cited as providing theoretical support for kernel methods. However, we nevertheless hope that the main ideas that make kernel methods attractive became clear. In particular, these include the fact that kernels address the following three major issues of learning and inference:

- They formalize the notion of *similarity* of data.
- They provide a *representation* of the data in an associated reproducing kernel Hilbert space.
- They characterize the *function class* used for estimation via the representer theorem [see equations (38) and (86)].

We have explained a number of approaches where kernels are useful. Many of them involve the substitution of kernels for dot products, thus turning a linear geometric algorithm into a nonlinear one. This way, one obtains SVMs from hyperplane classifiers, and kernel PCA from linear PCA. There is, however, a more recent method of constructing kernel algorithms, where the starting point is not a linear algorithm, but a linear criterion [e.g., that two random variables have zero covariance, or that the means of two samples are identical], which can be turned into a condition involving an efficient optimization over a large function class using kernels, thus yielding tests for independence of random variables, or tests for solving the two-sample problem. We believe that these works, as well as the increasing amount of work on the use of kernel methods for structured data, illustrate that we can expect significant further progress in the years to come.

**Acknowledgments.** We thank Florian Steinke, Matthias Hein, Jakob Macke, Conrad Sanderson, Tilmann Gneiting and Holger Wendland for comments and suggestions.

The major part of this paper was written at the Mathematisches Forschungsinstitut Oberwolfach, whose support is gratefully acknowledged. National ICT Australia is funded through the Australian Government's *Backing Australia's Ability* initiative, in part through the Australian Research Council.

T. Hofmann
Darmstadt University of Technology
Department of Computer Science
Darmstadt
Germany
E-mail: hofmann@int.tu-darmstadt.de

B. Schölkopf
Max Planck Institute
  for Biological Cybernetics
Tübingen
Germany
E-mail: bs@tuebingen.mpg.de

A. J. Smola
Statistical Machine Learning Program
National ICT Australia
Canberra
Australia
E-mail: Alex.Smola@nicta.com.au